\documentclass[11pt]{amsart}

\usepackage{amsthm,amsmath,stmaryrd,bbm,geometry,color}
\usepackage{amssymb}
\usepackage[english]{babel}
\usepackage[utf8]{inputenc}
\usepackage{graphicx}
\usepackage{verbatim}
\usepackage{enumitem}
\usepackage{mathtools}
\usepackage[titletoc]{appendix}
\usepackage{bm}

\usepackage[dvipsnames]{xcolor}
\usepackage[hyperindex=true,frenchlinks=true,colorlinks=true,
citecolor=Mahogany,linkcolor=Red,urlcolor=Tan,linktocpage]{hyperref}

\usepackage{tikz}
\usetikzlibrary{shapes}
\usetikzlibrary{fit}
\usetikzlibrary{decorations.pathmorphing}

\title{On the weak invariance principle for strictly stationary sequences in some Besov spaces}

\author{Davide Giraudo}
\address{Institut de Recherche Mathématique Avancée,
UMR 7501, Université de Strasbourg and CNRS,
7 rue René Descartes,
67000 Strasbourg, France}
\email{dgiraudo@unistra.fr}

\author{Sadillo Sharipov}
\address{V.~I. Romanovskiy Institute of Mathematics,
Uzbekistan Academy of Sciences,
9 Universitet Str.,
Tashkent, Uzbekistan}
\email{sadi.sharipov@yahoo.com}

\thanks{\textit{Funding.} This work was partially 
supported by a grant from the International Mathematical Union's Commission for 
Developing Countries (IMU-CDC) and the Simons Foundation.}

\numberwithin{equation}{section}
\renewcommand{\leq}{\leqslant}
\renewcommand{\geq}{\geqslant}

\newtheorem{Theorem}{Theorem}[section]
\newtheorem{Proposition}[Theorem]{Proposition}
\newtheorem{Lemma}[Theorem]{Lemma}
\newtheorem{Definition}[Theorem]{Definition}
\newtheorem{Corollary}[Theorem]{Corollary}

\theoremstyle{remark}
\newtheorem{Remark}[Theorem]{Remark}

\tikzstyle{Vertex}=[circle,draw=LimeGreen!80,fill=LimeGreen!8,
inner sep=1pt,minimum size=2mm,line width=1pt,font=\scriptsize]
\tikzstyle{Node}=[Vertex,draw=RoyalBlue!80,fill=RoyalBlue!8,inner sep=1.5pt]
\tikzstyle{Leaf}=[rectangle,draw=Black!70,fill=Black!16,
inner sep=0pt,minimum size=1mm,line width=1.25pt]
\tikzstyle{Edge}=[Maroon!80,cap=round,line width=1pt]
\tikzstyle{Mark1}=[draw=BrickRed!80,fill=BrickRed!8]
\tikzstyle{Mark2}=[draw=BurntOrange!80,fill=BurntOrange!8]
\tikzstyle{EdgeRew}=[->,RedOrange!80,cap=round,thick]

\newcommand{\Aca}{\mathcal{A}}

\newcommand{\Fca}{\mathcal{F}}

\newcommand{\Ica}{\mathcal{I}}

\newcommand \ens[1]{\left\{ #1\right\}}

\newcommand \R{\mathbb R}

\newcommand \PP{\mathbb P}
\newcommand{\el}{\mathbb L}
\newcommand{\E}[1]{\mathbb E\left[#1\right]}

\newcommand \Z{\mathbb Z}

\newcommand \abs[1]{\left|#1\right|}
\newcommand \eps{\varepsilon}

\newcommand{\pr}[1]{\left(#1\right)}
\newcommand{\prbig}[1]{\Big(#1\Big)}
\newcommand{\norm}[1]{\left\lVert #1 \right\rVert}
\newcommand{\normbig}[1]{\big\lVert #1 \big\rVert}

\newcommand{\ent}[1]{\left\lfloor #1\right\rfloor}

\begin{document}

\maketitle

\begin{abstract}
In this paper, we investigate the weak invariance principle in some Besov spaces for strictly stationary sequences. In particular, we provide sufficient  conditions on the moments and the dependence of a strictly stationary sequence to guaranty the weak invariance principle in these Besov spaces. Several dependence structures are considered, including projective conditions, $\alpha$-dependent, $\rho$-mixing, and $\tau$-dependent conditions. The main tools in the proofs are martingale approximation and deviation inequalities for weakly dependent sequences. Our results cover a large range of strictly stationary sequences that were not treated previously in the literature. Moreover, the presented conditions are sharp.
\end{abstract}

\section{Introduction}
Let $(\Omega,\mathcal{A},\mathbb{P})$ be a probability space, and 
$\mathbb{T}:\Omega\to\Omega$ be a bijective bi-measurable transformation preserving
the probability $\mathbb{P}$. Let $\mathcal{F}_0$ be a sub-$\sigma$-algebra
of $\mathcal{A}$ satisfying $\mathcal{F}_0\subseteq \mathbb{T}^{-1}(\mathcal{F}_0)$  then we define a non-decreasing filtration $\pr{\Fca_i}_{i\in\Z}$ by $\Fca_i= \mathbb{T}^{-i}\Fca_0$. We denote by $\Ica$ the $\sigma$-algebra of the $\mathbb{T}$-invariant sets, that is, the sets $A\in\Aca$ such that $\mathbb{T}^{-1}A=A$.
Given a random variable $X_0\colon\Omega\to \R$, we define the strictly stationary sequence
$\mathbb{X}=(X_i)_{i\in\Z}$ by $X_i=X_0\circ \mathbb{T}^i$. For any real $q\geq 1$, denote by $\el_q$ the  space of $\R$-valued random variables $Y$ such that
$\norm{Y}_q:=\big(\E{\abs{Y}^q}\big)^{1/q}$ is finite. 

Define the partial sums 
$$
S_{0,X_{0}}:=0, \ \ S_{k,X_{0}}:=\sum_{i=1}^{k}X_{i}, \ \ k\geq 1.
$$

Consider the random polygonal line with vertices $\left(k/n,n^{-1/2} S_{k,X_{0}}\right)$, $k=0,1,\ldots,n$:
\begin{equation}\label{eq:def_partial_sum_process}
W_{n,\mathbb{X}}(t)
:=\frac 1{\sqrt{n}}
S_{\lfloor nt \rfloor,X_{0}}
+\frac{
\bigl(nt-\lfloor nt\rfloor\bigr)}{\sqrt{n}}
X_{{\lfloor nt  \rfloor}+1}, \ \ t\in [0,1],
\end{equation}
where $\ent{x}$ denotes the lower part of $x\in\R$.

Donsker \cite{MR40613} has shown that when the sequence $\mathbb{X}$ is
independent and identically distributed (i.i.d.), centered and has variance
one, then the process
$\left(W_{n,\mathbb{X}}(t),\ t\in[0,1]\right)_{n\geq 1}$
converges in distribution in the space of continuous function $C([0,1])$, endowed with the uniform norm (denoted in the sequel by $\norm{\cdot}_\infty$), to a
standard Wiener process. Since then, considerable attention has been devoted to extending this result to dependent classes of random variables and to function spaces equipped with stronger topologies. We refer the readers to the monograph by Merlevède, Peligrad and Utev \cite{zbMATH07035947} for recent results on the weak invariance principles (WIP) for dependent random variables.

It should be noted that the weak convergence in $C([0,1])$ ensures convergence with respect to the uniform norm, which, in general, does not provide sufficient control over the local regularity of sample paths. The natural question then arises of whether the weak convergence of a properly normalized process \eqref{eq:def_partial_sum_process} remains valid in function spaces equipped with stronger topologies. This question motivates us to study the WIP in some H\"older and Besov spaces. We now introduce H\"older spaces.

For $\gamma\in (0,1]$, we consider the H\"olderian modulus of smoothness
of a function $x:[0,1]\to\mathbb{R}$ defined by
$$
\omega_{\gamma}(x,\delta)
:=
\sup_{\substack{s,t\in[0,1]\\0<t-s<\delta}}
\frac{|x(t)-x(s)|}{|t-s|^{\gamma}},
\qquad 0<\delta\leq 1.
$$
The H\"older space $\mathcal{H}_{\gamma}^{o}([0,1])$ is defined as
$$
\mathcal{H}_{\gamma}^{o}([0,1])
:=
\Big\{
x\in C([0,1]):
\lim_{\delta\to 0}\omega_{\gamma}(x,\delta)=0
\Big\}.
$$
$\mathcal{H}_{\gamma}^{o}([0,1])$ is a separable Banach space when endowed with the norm
$$
\|x\|_{\mathcal{H}_{\gamma}^{o}([0,1])}
:=
|x(0)|+\omega_{\gamma}(x,1).
$$
Invariance principles in Hölder spaces have attracted considerable attention in the literature. Early contribution in this field was due to Lamperti \cite{MR0143245}, who established the WIP in $\mathcal{H}_{\gamma}^{o}([0,1]), \gamma\in\pr{0,1/2}$. An alternative proof of this result was later given by Kerkyacharian and Roynette \cite{zbMATH00029174}, using a different approach based on Ciesielski's \cite{MR132389} characterization of H\"older spaces in terms of triangular functions. We stress that systematic investigations on limit theorems in $\mathcal{H}_{\gamma}^{o}([0,1]), \gamma\in\pr{0,1/2} $, were conducted by Hamadouche \cite{zbMATH01484659}, and Ra\v{c}kauskas and Suquet \cite{zbMATH02108226,zbMATH02118866,zbMATH02117810,zbMATH02083467,zbMATH07201347} and see the references therein. Mixing condition as well as $\tau$-dependence conditions were addressed in \cite{MR3615086}, while martingale approximation has been investigated in \cite{MR3426520,zbMATH06973953}. Dedecker and Merlevède \cite{MR3583992} established the WIP in $\mathcal{H}_{\gamma}^{o}([0,1])$ under an assumption involving the $\alpha$-dependence coefficients
sequence of $\mathbb{X}$, with applications to intermittent maps. In particular, this improves one of the result from \cite{MR3615086} valid under an assumption on the $\alpha$-mixing coefficients in the sense of Rosenblatt.
Alouin and Bigot \cite{zbMATH08101846} proved a Hölderian WIP for nonuniformly expanding maps, including bounded variation observables for certain intermittent maps.

Compared with H\"older spaces, Besov spaces provide a more flexible framework, due to the fact that they measure the regularity of sample paths through $\mathbb{L}_p$-control of their increments rather than uniform control. Let us define the Besov space on which we will work in the sequel.

For $p\geq 1$, define $\el_p\pr{[0,1]}$ as the space of (equivalence classes for the almost everywhere equality) of functions $x\colon [0,1]\to\R$ such that
$\|x \|_{\el_p\pr{[0,1]}}:= \pr{\int_{0}^{1}\abs{x\pr{t}}^{p}\mathrm{d}t}^{1/p}<\infty$.
The modulus of regularity of a function $x\in\el_{p} \pr{[0,1]}$ is defined by
\begin{equation*}
\omega_{p} \pr{x,\delta}:=\sup_{\abs{h}\leq \delta}\prbig{\int_{I_h}\abs{x\pr{u+h}-x\pr{u}}^{p} \mathrm{d}u}^{1/p}, \ \ 0<\delta\leq 1,
\end{equation*}
where for $h\in[-1,1]$, $I_h=[0,1]\cap [-h,1-h]$. For $\gamma\in (0,1]$ and $p\geq 1$, the Besov
space
$B_{p,\gamma}^{o}([0,1])$  is defined by
\begin{equation*}
B_{p,\gamma}^{o}([0,1]):=\Big\{x\in\el_p\pr{[0,1]}:\lim_{\delta\to 0} \delta^{-\gamma}\omega_p \pr{x,\delta}=0 \Big\},
\end{equation*}
and is endowed with the norm
\begin{equation*}
\norm{x}_{B_{p,\gamma}^{o}([0,1])}:=\|x\|_{\el_p\pr{[0,1]}}+\sup_{0<\delta<1}\delta^{-\gamma}\omega_p\pr{x,\delta}.
\end{equation*}
In this way, $B_{p,\gamma}^{o}([0,1])$ is a separable Banach space. Moreover, for
$0<\beta<\gamma$, the embedding  $B_{p,\gamma}^{o}([0,1]) \hookrightarrow B_{p,\beta}^{o}([0,1])$ is continuous. For $1\leq q<p$, the embedding  $B_{p,\gamma}^{o}([0,1]) \hookrightarrow
B_{q,\gamma}^{o}([0,1])$ is also continuous. Besov spaces were first introduced by O.~V.~Besov in \cite{zbMATH03158244}, who defined a family of function spaces in terms of moduli of smoothness and approximation properties via embedding and extension theorems. In the present paper, we work with the space $B_{p,\gamma}^{o}([0,1])$ in the ranges $p\gamma >1$, $
 \gamma\in\pr{1/p,1/2} $ and also $p\gamma\leq 1, \gamma\in\pr{0,1/2}$. These Besov spaces are  different according to the cases $p\gamma >1$ or $p\gamma\leq 1$.
Note that the restriction $\gamma>1/p$ ensures sufficient path regularity and it is related to the embedding of the corresponding Besov space into $C\pr{[0,1]}$, while restriction $\gamma<1/2$ is consistent with the regularity of Wiener process sample paths. It is worth to indicating that when $p\gamma >1$, an important advantage of the Besov spaces is that they admit a sequential characterization in terms of dyadic oscillations, or equivalently Faber-Schauder coefficients (see Theorem~\ref{thm:Schauder_decomposition}). When $p\gamma\leq 1$, different methods are required since the corresponding Besov space does not admit a Schauder decomposition.
We refer the reader to \cite{MR2319598} for more information about Besov spaces.

To the best of the authors' knowledge, the study of WIP in Besov spaces has attracted considerably less attention compared to Hölder spaces. In particular, Morel \cite{zbMATH02106939} proved a WIP in another sub-class of Besov spaces for sequences of i.i.d.\ random variables.
Giraudo and Račkauskas \cite{MR3736194} established WIP in $B_{p,\gamma}^{o}([0,1])$ for polygonal partial sum processes generated by stationary martingale difference sequences and as well as i.i.d.\ case. Their approach exploits a sequential characterization of Besov spaces in terms of Faber-Schauder coefficients and reduces tightness in the Besov topology to suitable estimates for dyadic blocks of partial sums. Let us be more specific 
about the results obtained in \cite{MR3736194}.
 In the regime $p\gamma > 1$, even in the i.i.d.\ case, finiteness of higher moments of the common distribution is required (see \cite{MR3736194}). More precisely, the WIP in $B_{p,\gamma}^{o}([0,1])$ holds if and only if
\begin{equation}\label{eq:tail_cond_iid}
    \lim_{t\to\infty}t^{q\pr{p,\gamma}}\PP\pr{\abs{X_0}\geq t}=0,
\end{equation}
where 
\begin{equation}\label{eq:def_qpalpha}
q\pr{p,\gamma}:=\prbig{\frac {1}{2}-\gamma+\frac{1}{p}}^{-1}.
\end{equation}
Notice that since $p\gamma > 1$ and $\gamma\in (0,1/2)$, one has $2<q\pr{p,\gamma}<p$.

 It turns out that for strictly stationary martingale difference sequences, condition \eqref{eq:tail_cond_iid} is not sufficient, but $X_0\in\el_{q\pr{p,\gamma}}$ is. Let us mention two recent results. Mies \cite{mies2026besovorliczmodulibrownianmotion} derived the WIP for moduli of regularity close to the critical, requiring sub-Gaussian bounds for the tails instead of \eqref{eq:tail_cond_iid}.
Giraudo and Sharipov \cite{giraudo2026weakinvarianceprincipletriangular} proved that polygonal line processes associated with partial sum processes of a triangular array of row-wise independent random variables converge in distribution in $B_{p,\gamma}^{o}([0,1])$ to a standard Wiener process. This naturally raises the question of finding tractable sufficient conditions for the validity of the WIP in $B_{p,\gamma}^{o}([0,1])$ in the strictly stationary setting.

The main purpose of the present paper is to establish the weak convergence of the process defined in \eqref{eq:def_partial_sum_process} in the space $B_{p,\gamma}^{o}([0,1])$, where
$\mathbb{X}$ is strictly stationary. More precisely, in the regime $p\gamma\leq 1$, we demonstrate that the WIP in $B_{p,\gamma}^{o}([0,1])$ reduces to the convergence of the
finite-dimensional distributions and uniform integrability of maxima of  the squares of normalized partial sums.
In the regime $p\gamma>1$, we provide sufficient conditions on the moments and the dependence of $\mathbb{X}$ involving $p$ and $\gamma$
which guaranty the WIP in $B_{p,\gamma}^{o}([0,1])$.
We show that summability of a series of norms of projectors of $X_0$ implies the WIP, extending the case of strictly stationary martingale difference sequences.
We also address the case of $\alpha$-dependent sequences, $\rho$-mixing sequences  as well as $\tau$-dependent sequences.
In the proofs, we mainly focus on two approaches: one is based on approximation
by martingales and the other one rests on probability inequalities
for weakly dependent sequences. The range of applications of our results covers a wide class of dependent stationary sequences: functions of an i.i.d.\ sequence like Volterra processes or functions of linear processes and some examples of the aforementioned weakly dependent sequences.

The paper is organized as follows. In Section~\ref{sec:main_results},  we give the main results, starting by the case $p\gamma\leq 1$, then we address the case $p\gamma>1$ by giving the Hannan type condition, then the results on the aforementioned dependence structure. In Section~\ref{sec:tools_Besov}, we collect the main tools to prove the WIP in $B_{p,\gamma}^o\pr{[0,1]}$, where $p\gamma >1$. Section~\ref{sec:proofs} is devoted to the proofs. Section~\ref{sec:deviation_inequalities} contains an appendix that collects the deviation inequalities for weakly dependent sequences that we will use in the sequel.
\section{Main Results}\label{sec:main_results}
\subsection{The case $p\gamma\leq 1$}

For this parameter range, as in the martingale difference setting considered in \cite{MR3736194}, no moments greater than the second are needed. In
order to guaranty tightness of the partial sum process, one needs to assume
uniform integrability of the squares of normalized partial sums, which holds
under various dependence conditions.
\begin{Theorem}\label{thm:Besov_alpha_p_leq_1}
Let $p$ and  $\gamma$ be such that  $p\geq 1$ and $p\gamma\leq 1$. Suppose that $\mathbb{X}=\pr{X_i}_{i\in \Z}$ is a strictly stationary centered sequence such that
the finite-dimensional distributions of $\pr{W_{n,\mathbb{X}}\pr{t},t\in [0,1]}_{n\geq 1}$ converge in distribution to
those of $\big(\widetilde{W}\pr{t},t\in [0,1]\big)$ and \begin{equation}\label{eq:uniform_integrability_for_p_alpha_leq_1}
\text{the sequence }\prbig{\frac 1n\max_{1\leq k\leq n}
    S_{k,X_0}  ^2}_{n\geq 1}\text{is uniformly
integrable}.
\end{equation}
Then $\pr{W_{n,\mathbb{X}}\pr{t},t\in [0,1]}_{n\geq 1}$ converges in distribution in
$B^{o}_{p,\gamma}\pr{[0,1]}$ to $\big(\widetilde{W}\pr{t},t\in [0,1]\big)$.
\end{Theorem}
\begin{Remark}
Note that the limiting process in this result is not necessarily a Wiener process. 
For example,  the limiting process may be of the form $\eta \cdot W\pr{t}$, where $\pr{W\pr{t},t\in [0,1]}$ is a standard Wiener process and $\eta$ is a $\mathbb{T}$-invariant map as in the case of non-ergodic martingale difference sequences (see \cite[Theorem~1]{MR3736194}).
\end{Remark}
Let us now give  some classes of strictly stationary sequences that
satisfy the assumptions of Theorem~\ref{thm:Besov_alpha_p_leq_1}. Under the next three conditions, the limiting process is of the form $\pr{\eta W\pr{t},t\in [0,1]}$, where $\pr{ W\pr{t},t\in [0,1]}$ is a standard Wiener process independent of $\eta$ and $\eta$ is $\mathbb{T}$-invariant.
 \begin{itemize}
     \item Hannan condition: suppose that $X_0$ is $\Fca_\infty$-measurable,  
$\E{X_0\mid \Fca_{-\infty}}=0$,  where $\Fca_{\infty}$ is the $\sigma$-algebra generated by $\bigcup_{i\in\Z} \Fca_i$ and 
$\Fca_{-\infty}=\bigcap_{i\in\Z} \Fca_i$
     and
     \begin{equation}\label{eq:Hannan_dim_1}
         \sum_{i\in\Z}
         \norm{\E{X_0\mid \Fca_{i}}-\E{X_0 \mid \Fca_{i-1}}}_2 <\infty.
     \end{equation}
     Then $\pr{W_{n,\mathbb{X}}\pr{t},t\in [0,1]}_{n\geq 1}$ converges in distribution in
$C\pr{[0,1]}$.
     Heyde \cite{zbMATH03466281} established a central limit theorem. When $X_{0}$ is
$\Fca_0$-measurable and $\mathbb{T}$ is weakly mixing, Hannan \cite{MR0331683} proved 
that $\pr{W_{n,\mathbb{X}}\pr{t},t\in [0,1]}_{n\geq 1}$ converges in distribution in
$C\pr{[0,1]}$. Dedecker and Merlevède \cite[Corollary~3]{MR2019054} showed that the weak mixing assumption on $\mathbb{T}$ can be dropped, while Dedecker, Merlevède and Voln\'{y} \cite{MR2359065} proved that $X_0$ may be taken to be $\Fca_{\infty}$-measurable instead of $\Fca_{0}$-measurable. 
     \item Maxwell and Woodroofe condition: it reads as 
     \begin{equation}\label{eq:MW_dim_1}
         \sum_{k=1}^{\infty}\frac 1{k^{3/2}}
         \prbig{\normbig{\sum_{i=0}^{k-1}\E{X_i\mid\Fca_0}  }_2
         +\normbig{\sum_{i=0}^{k-1}\pr{ X_{-i}-\E{X_{-i}\mid\Fca_0}
}}_2}<\infty.
     \end{equation}
     Then $\pr{W_{n,\mathbb{X}}\pr{t},t\in [0,1]}_{n\geq 1}$ converges in distribution in
$C\pr{[0,1]}$. The central limit theorem was derived by Maxwell and Woodroofe 
\cite{MR1782272} when $X_0$ is $\Fca_0$-measurable (note that in this case, the
second term in \eqref{eq:MW_dim_1} vanishes).  Peligrad and Utev \cite{MR2123210} obtained
     its functional form using an appropriate maximal
inequality, but still under the assumption that $X_0$ is $\Fca_0$-measurable, which 
has been removed by Voln\'{y} \cite{MR2344817}.
     \item Dedecker and Rio  projective condition: if
     $X_0$ is $\Fca_0$-measurable and the sequence\\
$\pr{X_0\E{S_{n,X_0}\mid\Fca_0}}_{n\geq 1}$ converges in $\el_1$, then
$\pr{W_{n,\mathbb{X}}\pr{t},t\in [0,1]}_{n\geq 1}$ converges in distribution in
$C\pr{[0,1]}$. This projective condition was introduced by Dedecker and Rio \cite{MR1743095}.
 \end{itemize}
         Under the previous projective conditions, assumption \eqref{eq:uniform_integrability_for_p_alpha_leq_1} is satisfied. Indeed, under \eqref{eq:Hannan_dim_1} and \eqref{eq:MW_dim_1}, it is known \cite{Gordin_Peligrad} that there exists a strictly stationary martingale difference sequence $\pr{D_i}_{i\in\Z}$ such that 
         \begin{equation*}
            \displaystyle\lim_{n\to\infty}
         n^{-1}\E{\max_{1\leq k\leq n}(S_{k,X_0}-S_{k,D_0} )^2}= 0.
         \end{equation*}
         It is standard that 
         \eqref{eq:uniform_integrability_for_p_alpha_leq_1} holds for strictly stationary martingale difference sequences.  Under the  Dedecker and Rio  projective condition, 
         \eqref{eq:uniform_integrability_for_p_alpha_leq_1} is a consequence of Proposition~1 in \cite{MR1743095}.
         
        Therefore, the weak convergence in $B_{p,\gamma}^o\pr{[0,1]}$ also holds in the case where $p\gamma\leq 1$, $ \gamma\in\pr{0,1/2}$.
 Note that all these projective conditions are trivially satisfied if
$\mathbb{X}$ is a martingale difference sequence with respect
to
 $\pr{\Fca_i}_{i\in\Z}$.
 It was shown by Durieu and Voln\'{y} \cite{zbMATH05547583,MR2446326} that these
criteria are independent: taking any pair of them, one can find $\mathbb{X}$ satisfying one of the criterion but not the other.
We refer the reader to the aforementioned paper and the references therein for explicit examples
satisfying the previous conditions.

\subsection{A Hannan type condition}

Our first result for the regime $p\gamma>1$ establishes a WIP under a condition analogous to \eqref{eq:Hannan_dim_1}, involving a series of norms of projections.
\begin{Theorem}\label{thm:Hannan}
Let $p>2$,  $\gamma\in\pr{1/p,1/2}$ and $q\pr{p,\gamma}$ be defined in \eqref{eq:def_qpalpha}. 
Suppose that $\mathbb{X}=\pr{X_i}_{i\in\Z}$ is a strictly stationary sequence such that $X_0$ is $\Fca_{\infty}$-measurable, $\E{X_0\mid\Fca_{-\infty}}=0$ and 
\begin{equation}\label{eq:Hannan}
\sum_{i\in\Z}\norm{\E{X_{0}\mid\Fca_{i}}-
\E{X_{0}\mid\Fca_{i-1}}  }_{q\pr{p,\gamma}}<\infty.
\end{equation}
Then $\pr{W_{n,\mathbb{X}}\pr{t},t\in [0,1]}_{n\geq 1}$ converges in distribution in
$B^o_{p,\gamma}([0,1])$ to $\pr{\eta W\pr{t},t\in [0,1]}$, where $\eta$ is measurable with respect to the $\sigma$-algebra of $\mathbb{T}$-invariant sets and is independent of the standard Wiener process $\pr{W\pr{t},t\in [0,1]}$.
\end{Theorem}
\begin{Remark}
    This result complements \cite[Theorem~2.6]{MR3426520} obtained in $\mathcal{H}_\gamma^o\pr{[0,1]}$, which corresponds to $B_{\infty,\alpha}^o\pr{[0,1]}$ (where the $\el_p$-norm is replaced by the $\el_\infty$-norm).
\end{Remark}

 Let us now give some applications of Theorem~\ref{thm:Hannan} to sequences of the form:
\begin{equation}\label{eq:def_Bernoulli_shift}
     X_i=f\pr{\pr{\eps_{i-k}}_{k\in\Z}}, \ \qquad {i\in \mathbb{Z}},
 \end{equation}
where $f\colon \R^{\Z}\to\R$ is measurable, where $\R^{\Z}$ is the set of $\R$-valued sequences, which is endowed with the product $\sigma$-algebra of the Borel subsets of $\R$
and $\pr{\eps_k}_{k\in\Z}$ is an i.i.d.\ sequence. Following Wu \cite{MR2172215}, we define the physical dependence coefficients of the sequence $\pr{X_i}$ by
\begin{equation}\label{eq:physical_dependence_measure}
\delta_{\mathbb{X},i,p}:=\norm{X_i-X_i^*}_p,
\qquad i\in\Z,\quad p\geq 1.
\end{equation}
Here, 
 $X_i^*=f\big(\big(\eps_{i-k}^*\big)_{k\in\Z}\big)$, $\eps^*_j=\eps_j$ if
$j\neq 0$ and $\eps^*_0=\eps'_0$, where $\eps'_0$ is independent of
$\pr{\eps_k}_{k\in\Z}$ and has the same distribution as $\eps_0$.

In the next corollary, we state that the conditions of Theorem~\ref{thm:Hannan} are satisfied under summability of the 
dependence coefficients. Note that sequences  of the form \eqref{eq:def_Bernoulli_shift} are necessarily ergodic, hence the limiting process is a Wiener process multiplied by a constant.
\begin{Corollary}\label{cor:Bernoulli}
Let $p>2$, $\gamma\in\pr{1/p,1/2}$ and $q\pr{p,\gamma}$ be defined in \eqref{eq:def_qpalpha}.
Let $\mathbb{X}=\pr{X_i}_{i\in\Z}$ be a centered sequence defined in \eqref{eq:def_Bernoulli_shift}. Suppose that the coefficients $\delta_{\mathbb{X},i,p}$
defined in \eqref{eq:physical_dependence_measure} satisfy
\begin{equation}\label{eq:cond_Bernoulli}
\sum_{i\in\Z}\delta_{\mathbb{X},i,q\pr{p,\gamma}}<\infty.
\end{equation}
Then $\pr{W_{n,\mathbb{X}}\pr{t},t\in [0,1]}_{n\geq 1}$ converges in distribution in
$B^o_{p,\gamma}([0,1])$ to $\pr{\sigma W\pr{t},t\in [0,1]}$,
where\\ $\pr{  W\pr{t},t\in [0,1]}$ is a standard Wiener process and 
\begin{equation}\label{eq:sigma_2_bernoulli}
    \sigma^2:=\sum_{k\in\Z}\operatorname{Cov}\pr{X_0,X_k}.
\end{equation}
\end{Corollary}
It is a standard fact that \eqref{eq:cond_Bernoulli} implies the absolute convergence of the series 
involved in \eqref{eq:sigma_2_bernoulli}.
Let us give three examples of sequences for which Corollary~\ref{cor:Bernoulli} can be applied (see \cite[pages 5966-5968]{MR3256190} for bounds on the coefficients $\delta_{\mathbb{X},i,q\pr{p,\gamma}}$).
\begin{enumerate}
    \item Hölder function of a linear process.
     Let $\mathbb{X}=\pr{X_i}_{i\in\Z}$ be a H\"{o}lder continuous function of a linear process as
considered, for example, in  \cite{MR1487431,MR2212691}.
More precisely, define
$X_i$ by
\begin{equation*}
X_i:=\varphi\big(\sum_{j\in\Z}a_j\varepsilon_{i-j}\big),\quad i\in\Z,
\end{equation*}
where $\varphi\colon\R\to\R$ is $\varrho$-H\"{o}lder continuous for some
$\varrho\in
(0,1]$, $\pr{a_j}_{j\in \Z}$ is a sequence of real numbers such that
$\displaystyle\sum_{j\in\Z} \abs{a_j}^{\varrho}<\infty$ and $\pr{\varepsilon_j}_{j\in\Z}$ is
an
i.i.d.\ sequence with $\varepsilon_0\in\el_{q\pr{p,\gamma}} $.
    \item Function of a Gaussian linear process. Assume   that $\mathbb{X}=\pr{X_{i}}_{i\in\Z}$ can be written as a function of a Gaussian linear process: let
$X_i=\varphi\pr{Y_i}$ with
 \begin{equation*}
Y_{i}:=\sum_{ j \in\Z  }a_{j}\varepsilon_{i-j},\quad i\in\Z,
\end{equation*}
 where $\varphi:\R\rightarrow\R$,
 $\pr{\varepsilon_{j}}_{j \in \Z}$ is an i.i.d.\ sequence,
$\varepsilon_{0}$ has
a standard normal distribution and $a_j \in \R$ for all $j\in\Z$ with $\displaystyle\sum_{ j
\in\Z
}\abs{a_{j}}<\infty$ as well as $\displaystyle\sum_{ j \in\Z
}a_{j}^2=1$.  Such processes were considered, e.g., in \cite{MR2498953}.
Given  that $\varphi\pr{Y_0}\in\mathbb L_2$
and
$\E{\varphi\pr{Y_{0}}}=0$, the following expansion holds:
\begin{equation*}
\varphi\pr{Y_{i}}=\sum_{q=1}^{\infty}c_q\pr{\varphi}H_q\pr{Y_{i}},
\end{equation*}
where the
$q$-th Hermite polynomial $H_q$ is defined by
\begin{equation*}
H_q\pr{x}:=\pr{-1}^q\exp\prbig{\frac{x^2}2}\frac{\mathrm{d}^q}{\mathrm{d}x^q}\exp\prbig{-\frac{x^2}2},
\end{equation*}
and
\begin{equation*}
c_q\pr{\varphi}:=\frac{1}{q!}\E{\varphi\pr{Y_{0}}H_q\pr{Y_{0}}}.
\end{equation*}
Since $\varphi\pr{Y_0}\in\mathbb L_2$, the series
$\displaystyle\sum_{q=1}^{\infty}q!\pr{c_q\pr{\varphi}}^2$ is finite.
 The conditions of Corollary~\ref{cor:Bernoulli} are satisfied if  
\begin{equation*}
\sum_{j\in\Z}\abs{a_j}<\infty
\text{ and }
 \sum_{q=1}^{\infty}\sqrt{q\cdot q!}
\abs{c_q\pr{\varphi }}<\infty.
\end{equation*}
    \item Volterra processes. Let $\mathbb{X}=\pr{X_{i}}_{i\in\Z}$ be a Volterra process 
    defined by 
\begin{equation*}
 X_{i}:=\sum_{j,j'\in\Z,j\neq j'}a_{j,j'}\varepsilon_{i-j}\varepsilon_{i-j'},\quad i\in\Z,
\end{equation*}
where $\pr{\varepsilon_j}_{j\in\Z}$ is an i.i.d.\ centered sequence with
$ \varepsilon_0\in\el_2$
and $\displaystyle\sum_{j,j'\in\Z,j\neq j'}a_{j,j'}^2<\infty$. This type of processes and their extensions
to Volterra series are considered in \cite{MR991969,MR713994}.

Under the assumptions 
\begin{equation*}
\varepsilon_0\in\el_{q\pr{p,\gamma}} \text{ and }   \sum_{j\in\Z}
\prbig{\sum_{j'\in\Z, j'\neq j} \pr{a_{j,j'}^2+a_{j',j}^2}}^{1/2}  <\infty,
\end{equation*}
 $\mathbb{X}$ satisfies the conditions of Corollary~\ref{cor:Bernoulli}.
\end{enumerate}

\subsection{$\alpha$-dependence condition}
The $\alpha$-dependence coefficients were introduced in
\cite{MR2682267}.
We follow the notations and definitions from 
\cite{MR2682267,MR3583992}.
\begin{Definition}\label{def:alpha_dep1}
For any integrable random variable $Z$, denote $Z^{\pr{0}}=Z-\E{Z}$. For any integrable random variable $V=\pr{V_1,\dots,V_\ell}$ with values
in $\R^\ell$ and any sub-$\sigma$-algebra $\Fca$ of $\Aca$, let
$$
\alpha\pr{\Fca,V}:=\sup_{x_1,\dots,x_\ell\in\R}
\Big\lVert \mathbb{E}\Big[\prod_{j=1}^\ell\pr{\mathbf{1}_{V_j\leq x_j }}^{\pr{0}} \mid\Fca \Big]
-\mathbb{E}\Big[\prod_{j=1}^\ell \pr{\mathbf{1}_{V_j\leq x_j }}^{\pr{0}} \Big]\Big\rVert_1.
$$
\end{Definition}
Given a strictly stationary sequence $\mathbb{Y}=\pr{Y_i}_{i\in\Z}$, we define the $\alpha$-dependence coefficients as 
$$
\alpha_{k,\mathbb{Y}}\pr{0}:=1/2 \mbox{ and }\alpha_{k,\mathbb{Y}}\pr{n}
:=\max_{1\leq \ell \leq k}\sup_{n\leq i_1\leq \dots\leq i_\ell}
\alpha\pr{\Fca_0,\pr{Y_{i_1},\dots,Y_{i_\ell}}},\quad  n\geq 1.
$$
In order to derive a WIP, one needs only assumptions on the sequences
$\pr{\alpha_{1,\mathbb{Y}}\pr{n}}_{n\geq 1}$ and $\pr{\alpha_{2,\mathbb{Y}}\pr{n}}_{n\geq 1}$
(see \cite{MR2682267,MR3583992}).
The stationary sequences in which we will be interested are some functions of $\alpha$-dependent sequences. Let us be more specific about the kind of functions.
\begin{Definition}\label{def:alpha_dep2}
We say that $H\colon \R_+ \to  [0,1]$
is a tail function if it is non-increasing, right-continuous, and converges to $0$ at infinity. The quantile function $Q=H^{-1}$ is the generalized inverse of the function $H$, that is, for $u\in [0,1]$, $Q\pr{u}=\inf\ens{t\geq 0: H\pr{t}\leq u}$ (with the convention $\inf \emptyset= \infty$).
\end{Definition}
\begin{Definition}\label{def:alpha_dep3}
Given a probability measure $\mu$ on $\R$ and $Q=H^{-1}$ an integrable quantile function, we define $\widetilde{\operatorname{Mon}}\pr{Q,\mu}$ as the set of the functions $g$ defined on $\R$ which are monotonic on some interval of $\R$
and null elsewhere, and such that for each $t\in\R_+$, $\mu\pr{\abs{g}>t}\leq H\pr{t}$. Finally, let $\widetilde{\Fca}\pr{Q,\mu}$ be the closure in $\el_1\pr{\mu}$ of the set of functions that can be
written as $\sum_{\ell=1}^L a_\ell f_\ell$, where $\sum_{\ell=1}^L \abs{a_\ell}\leq 1$ and $f_\ell\in \widetilde{\operatorname{Mon}}\pr{Q,\mu}$.
\end{Definition}

Our result for $\alpha$-dependent sequences reads as follows.
\begin{Theorem}\label{thm:alpha_dep}
Let $p>2$, $\gamma\in \pr{1/p,1/2+1/p-2/\pr{p+2}}$ and $q\pr{p,\gamma}$ be defined in \eqref{eq:def_qpalpha}. Let $Q$ be a quantile function in $\el_{q\pr{p,\gamma}}\pr{[0,1]}$,   $\mathbb{Y}=\pr{Y_i}_{i\in\Z}$ and  $\PP_{Y_0}$ denote the distribution of $Y_0$.
Assume $\mathbb{X}=\pr{X_i}_{i\in\Z}$ has the form
$X_i=f\pr{Y_i}-\E{f\pr{Y_i}}$, where $f\in\widetilde{\Fca}\pr{Q,\PP_{Y_0}}$. Define $R_{\pr{\alpha}} \pr{u}:=\alpha_{2,\mathbb{Y}}^{-1}\pr{u}Q\pr{u}$, $u\in [0,1]$. Suppose that
\begin{equation}\label{eq:cond_alpha_dep}
    \lim_{x\to\infty }x^{q\pr{p,\gamma}-1}\int_0^1Q\pr{u}\mathbf{1}_{R_{\pr{\alpha}} \pr{u}>x}\mathrm{d}u=0.
\end{equation}
Then $\pr{W_{n,\mathbb{X}}\pr{t},t\in [0,1]}_{n\geq 1}$ converges in distribution in
$B^{o}_{p,\gamma}([0,1])$ to $\pr{\sigma W\pr{t},t\in [0,1]}$,
where \\
$\pr{W\pr{t},t\in [0,1]}$ is a standard Wiener process, and $\sigma^2$ is given by \eqref{eq:sigma_2_bernoulli}.
\end{Theorem}
\begin{Remark}
 Here, the term $\sigma^2$ is given by a series of covariances, whose absolute convergence is established in  \cite[Theorem~1.13]{MR2682267}.   
\end{Remark}

\begin{Remark}\label{rem:recover_iid}
It will be seen in the proof that condition \eqref{eq:cond_alpha_dep} implies 
\eqref{eq:conv_for_tightness} which in turn implies \eqref{eq:tail_cond_iid}. When we apply Theorem~\ref{thm:alpha_dep} to an i.i.d.\ sequence, we may recover  \cite[Theorem~2]{MR3736194} when $p>2$, $\gamma\in \pr{1/p,1/2+1/p-2/\pr{p+2}}$. 
\end{Remark}
\begin{Remark}
 Notice that the regime $p\gamma >1$ is partially covered, since we have the restriction $\gamma <1/2+1/p-2/\pr{p+2}$. This comes from the application of Proposition~\ref{prop:deviation_alpha_dep}: one needs, as we shall see in the proofs for the other dependence structures, that the free parameter in the deviation inequality must be taken between $q\pr{p,\gamma}$ and $p$. However, in order to control the contribution of the last term of \eqref{eq:deviation_alpha_dep}, one needs also that $r/2 >q\pr{p,\alpha}-1$ and the inequality $p >2q\pr{p,\alpha}-2$ only holds in the aforementioned range.    
\end{Remark}
\begin{Remark}
    A bound for the $\alpha$-dependence coefficients for intermittent maps has been given in \cite{MR3583992}. Following the presented approach, one can also derive 
   a WIP in $B_{p,\gamma}^o\pr{[0,1]}$ for intermittent maps by letting the role played by $\pr{1/2-\gamma}^{-1}$ in  Hölder spaces by $q\pr{p,\gamma}$, provided that $p,\gamma$ satisfy the assumptions of Theorem~\ref{thm:alpha_dep}.
\end{Remark}
 
\subsection{$\rho$-mixing condition}

$\rho$-mixing coefficients were introduced by Hirschfeld in \cite{zbMATH03019401} (see also Kolmogorov and Rozanov \cite{zbMATH03149369}) and are defined by
\begin{equation*}
\rho_{\mathbb{X}}\pr{n}:=
\sup\ens{\frac{\abs{\operatorname{Cov}\pr{Y,Z}}}{\sqrt{\operatorname{Var}\pr{Y}\operatorname{Var}\pr{Z}  } }, Y\in \el_2\pr{\Fca_{-\infty}^0} , Z\in \el_2\pr{\Fca_{n}^{\infty}} , Y,Z\text{ are non-constant}    },
\end{equation*}
where $\Fca_{-\infty}^{0}=\sigma\pr{X_{i},i\leq 0}$ and $\Fca_{n}^{\infty}=\sigma\pr{X_{i}, i\geq n}$.

The next result requires the same moment condition as in the i.i.d.\ case, together with a standard decay condition on the $\rho$-mixing coefficients of $\mathbb{X}$.
\begin{Theorem}\label{thm:rho_mixing}
Let $p>2$, $\gamma\in\pr{1/p,1/2}$ and $
q\pr{p,\gamma}$ be defined in \eqref{eq:def_qpalpha}. Suppose that the strictly stationary centered  sequence $\mathbb{X}=\pr{X_i}_{i\in\Z}$ satisfies 
\eqref{eq:tail_cond_iid} and  
\begin{equation}\label{eq:series_coeff_rho_mixing}
 \sum_{k=0}^\infty \rho_{\mathbb{X}}\big(2^k\big)<\infty.
\end{equation}
Then $\pr{W_{n,\mathbb{X}}\pr{t},t\in [0,1]}_{n\geq 1}$ converges in distribution in
$B^o_{p,\gamma}([0,1])$ to $\pr{c W\pr{t},t\in [0,1]}$, where $\displaystyle c^2:=\lim_{n\to\infty}n^{-1}\E{S_{n,X_0}^2} \in [0,\infty)$  
and $\pr{W\pr{t},t\in [0,1]}$ is a standard Wiener process.
\end{Theorem}
\begin{Remark}
The existence of the limit defining $c^2$ follows from \cite[Theorem~4.1]{MR672297}.
\end{Remark}
\begin{Remark}
When we apply Theorem~\ref{thm:rho_mixing} to an i.i.d.\ sequence, we may recover \cite[Theorem~2]{MR3736194}.
\end{Remark}
  
\subsection{$\tau$-dependence condition}
In order to define the $\tau$-dependence coefficients of a stationary sequence, we first need a result about conditional probability (see \cite[Theorem 33.3]{MR1324786}).

\begin{Lemma}
Let $(\Omega, \Aca, \mathbb{P})$ be a probability space, $\mathcal{M}$ a sub-$\sigma$-algebra of $\mathcal{A}$, and $X$ an $\R$-valued random variable with distribution $\mathbb{P}_X$. Then there exists a measurable function $\mathbb{P}_{X|\mathcal{M}}\colon \mathcal{B}(\mathbb{R}) \times \Omega\to[0, 1]$ such that
\begin{enumerate}
    \item For any $\omega \in \Omega$, $\mathbb{P}_{X|\mathcal{M}}(\cdot, \omega)$ is a probability measure on $\mathcal{B}(\mathbb{R})$.
    \item For any $A \in \mathcal{B}(\mathbb{R})$, $\mathbb{P}_{X|\mathcal{M}}(A, \cdot)$ is a version of $\E{\mathbf{1}_{ X \in A } \mid \mathcal{M}}$.
\end{enumerate}
\end{Lemma}
We now recall the $\tau$-dependence coefficients, which were introduced by Dedecker and Prieur \cite{MR2199291}. We denote by $\Lambda_{1}(\mathbb{R})$ the set of 1-Lipschitz functions from $\mathbb{R}$ to $\mathbb{R}$ and define the quantity
$$
\mathcal{W}(\mathbb{P}_{X|\mathcal{M}}):= \sup \Big\{ \Big\lvert \int_{\R} f(x)\PP_{X|\mathcal{M}}(\mathrm{d}x)-\int_{\R} f(x)\mathbb{P}_X(\mathrm{d}x) \Big\rvert, f \in \Lambda_{1}(\mathbb{R}) \Big\}.
$$
For an integrable random variable $X$ and a sub-$\sigma$-algebra $\mathcal{M}$, the coefficient $\tau$ is defined by
$$
\tau(\mathcal{M}, X):=\Vert \mathcal{W}(\mathbb{P}_{X|\mathcal{M}}) \Vert_{1}.
$$
This definition can be extended to random variables with values in finite-dimensional vector spaces. If $\ell$ is a positive integer, we endow $\mathbb{R}^\ell$ with the norm $\Vert x-y \Vert:=\sum_{j=1}^{\ell} |x_{j}-y_{j}|$ and define $\Lambda_{1}(\mathbb{R}^\ell)$ as the set of 1-Lipschitz functions from $\mathbb{R}^\ell$ to $\mathbb{R}$.

\begin{Definition}
Let $(\Omega, \Aca, \mathbb{P})$ be a probability space, $\mathcal{M}$ be a sub-$\sigma$-algebra of $\Aca$, and $X$ be an $\mathbb{R}^\ell$-valued random variable. Define
$$
\tau(\mathcal{M}, X):= \sup \Big\{\tau(\mathcal{M}, f(X)): f \in \Lambda_1\big(\R^\ell\big) \Big\}.
$$
\end{Definition}

Introduce the $\tau$-dependence coefficients for a sequence of $\R$-valued random variables.

\begin{Definition}\label{def1.4}
Let $\mathbb{X}=\pr{X_i}_{i\in\Z}$ be a sequence of random variables and let $\mathcal{M}_{p}:=\sigma\pr{X_{k}, k \leq p}$. For any integer $i\geq 1$, define the $\tau$-dependence coefficient $\tau_{\mathbb{X}}\pr{i}$ by 
$$
\tau_{\mathbb{X}}\pr{i}:=\sup_{p,\ell\geq 1} \frac{1}{\ell} \sup \left\{\tau(\mathcal{M}_{p}, (X_{j_{1}}, \ldots, X_{j_\ell})): p+i \leq j_{1}< \cdots <j_{\ell} \right\}.
$$
\end{Definition} 
For an $\R$-valued random variable $X$, we denote by  $Q_X(\cdot)$ the inverse function defined by $Q_{X}(u):=\inf \ens{t>0: \PP\pr{\abs{X}>t} \leq u}$. If $(\delta_i)_{i\geq 0}$ is a non-increasing sequence of non-negative numbers, we define $\delta^{-1}(u):=\inf \{k \geq 0: \delta_{k} \leq u\}$. Let $G_{X_0}$ be the inverse of map $x\mapsto \int_0^xQ_{X_0}\pr{u}\mathrm{d}u$ and $R_{\pr{\tau}} := ((\tau_{\mathbb{X}}/2)^{-1} \circ G_{X_{0}}^{-1}) Q_{X_{0}}$. 

We are now in position to state our result on $\tau$-dependent sequences.
\begin{Theorem}\label{thm:tau_dependent}
Let $p>2$, $\gamma\in\pr{1/p,1/2}$, and $q\pr{p,\gamma}$ be defined in \eqref{eq:def_qpalpha}. Let $\mathbb{X}=(X_{i})_{i\in \mathbb{Z}}$ be a strictly stationary centered sequence such that
\begin{equation}\label{eq:main-condition for tau}
\lim_{t\to\infty} t^{q(p,\gamma)-1} \int_{0}^{1} Q_{X_{0}}(v) \mathbf{1}_{R_{\pr{\tau}}(v)>t  } \, \mathrm{d}v=0.
\end{equation}
Then $\pr{W_{n,\mathbb{X}}\pr{t},t\in [0,1]}_{n\geq 1}$ converges in distribution in
$B^{o}_{p,\gamma}([0,1])$ to $\pr{\sigma W\pr{t},t\in [0,1]}$, where $  \sigma^{2}$ is given by \eqref{eq:sigma_2_bernoulli} and $\pr{W\pr{t},t\in [0,1]}$ is a standard Wiener process.
\end{Theorem}
\begin{Remark}
  Notice that condition \eqref{eq:main-condition for tau} implies finiteness of 
$\int_0^1 R_{\pr{\tau}}\pr{u}Q_{X_0}\pr{u}\mathrm{d}u$, which in turn  implies condition (5.5) from \cite{MR2105738} and guarantees the convergence of the series defining $\sigma^2$.   
\end{Remark}
\begin{Remark}
When we apply Theorem~\ref{thm:tau_dependent} to an i.i.d.\ sequence, we may recover \cite[Theorem~2]{MR3736194}.
\end{Remark}

\section{Tools for Besov spaces}\label{sec:tools_Besov}
We present some tools in order to handle the convergence in distribution of processes in Besov spaces.
First, we introduce the notation for the Faber--Schauder functions. A convenient way to parametrize this family is by using dyadic numbers directly.
We denote by $\mathbf{D}_{j}$ the set of dyadic numbers in $[0,1]$ of level $j$, that is,
$$
\mathbf{D}_{0}=\{0,1\}, \qquad
\mathbf{D}_{j} =\left\{(2\ell-1)2^{-j} \; ; \; 1 \leq \ell \leq 2^{j-1} \right\}, \quad j 
\geq 1.
$$
Set
$$
\mathbf{D} =\bigcup_{j \geq 0}\mathbf{D}_{j}
$$
and, for $r \in \mathbf{D}_j$, write
$$
r^{-}:=r-2^{-j}, \qquad r^{+}:=r+2^{-j}.
$$
The triangular Faber--Schauder functions $\Lambda_r$ for $r \in \mathbf{D}_j$, $j>0$, are defined by
$$
\Lambda_r(t)=
\begin{cases}
2^j(t-r^{-}), & t \in (r^{-},r],\\[1mm]
2^j(r^{+}-t), & t \in (r,r^{+}],\\[1mm]
0, & \text{otherwise}.
\end{cases}
$$
When $j=0$, we just take the restriction to $[0,1]$:
$$
\Lambda_{0}(t)=1-t, \qquad \Lambda_{1}(t)=t, \quad t\in[0,1].
$$
\begin{Theorem}[Ciesielski et al \cite{MR132389}]
\label{thm:Schauder_decomposition}
Let $p>1$ and $\gamma\in\pr{1/p,1}$. The Faber--Schauder system $(\Lambda_{r})_{r\in \mathbf{D}}$ is a Schauder basis
for $B^{o}_{p,\gamma}([0,1])$: each $x\in B^{o}_{p,\gamma}([0,1])$ has a unique representation
$$
x=\sum_{r\in\mathbf{D}} \lambda_{r}(x)\Lambda_{r},
$$
where
$$
\lambda_{r}(x):=\frac{2x(r)-x(r^{+})+x(r^{-})}{2}, \qquad r\in \mathbf{D}_{j},\ j\geq 1,
$$
and, in the particular case $j=0$,
$$
\lambda_{0}(x):=x(0), \qquad \lambda_{1}(x):=x(1).
$$
Moreover, the norm is equivalent to the sequential norm:
\begin{equation}\label{sequential norm}
\|x\|_{B_{p,\gamma}^{o}([0,1])} \sim \|x\|_{B_{p,\gamma}^{o}([0,1])}^{\mathrm{seq}}
:= \sup_{j\geq 0} 2^{j(\gamma-1/p)}
\prbig{\sum_{r\in \mathbf{D}_j} |\lambda_{r}(x)|^{p}}^{1/p}.
\end{equation}
\end{Theorem}

The Schmidt orthogonalization procedure (with respect to the inner product in $\mathbb{L}_{2}([0,1])$) applied to the
Faber--Schauder system leads to the Franklin system $(f_{k})_{ k\geq 0}$:
$$
f_{k}(t)=\sum_{i=0}^{k} c_{ik}\Lambda_{r_{i}}(t), \qquad t\in [0,1],
$$
with $c_{kk}>0$ for $k\geq 0$, where the matrix $(c_{ik})$ is uniquely determined and
$(r_{i})_{i\geq 0}$ is an enumeration of $\mathbf{D}$.

Let us conclude this section by stating a tightness criterion in $B_{p,\gamma}^o\pr{[0,1]}$, derived by following the approach presented in \cite{MR1736910}.
 \begin{Proposition}[{\cite[Theorem~10]{MR3736194}}]\label{prop:tightness_criterion} Let $p$ and $\gamma$ be such that
 $\gamma\in\pr{1/p,1/2}$. Suppose that $\mathbb{X}=\pr{X_i}_{i\in\Z}$ is a
 strictly stationary sequence which satisfies
 \begin{equation}\label{eq:conv_for_tightness}
 \lim_{J\to\infty}\limsup_{N\to\infty}\sum_{j=J}^N
 2^j\int_0^1 u^{p-1}\PP\prbig{\max_{1\leq k\leq 2^{N-j}}\abs{S_{k,X_0}}>
 2^{N/2}2^{j(1/p -\gamma)}  u }\mathrm{d}u=0.
 \end{equation}
 Then the sequence $\pr{W_{n,\mathbb{X}}\pr{\cdot}}_{n\geq 1}$
 is tight in $B^o_{p,\gamma}([0,1])$.
 \end{Proposition}
It is worth noticing that \eqref{eq:conv_for_tightness} implies \eqref{eq:tail_cond_iid}. Indeed,
\begin{align*}
0&= \lim_{J\to\infty}\limsup_{N\to\infty}\sum_{j=J}^N
 2^j\int_0^1 u^{p-1}\PP\prbig{\max_{1\leq k\leq 2^{N-j}}\abs{S_{k,X_0}}>
 2^{N/2}2^{j(1/p -\gamma)}  u }\mathrm{d}u\\
 &\geq \limsup_{N\to\infty}2^N\int_{0}^1 u^{p-1}\PP\prbig{ \abs{X_0}>
 2^{N/q\pr{ p,\gamma} } u }\mathrm{d}u \\
 &\geq \limsup_{N\to\infty}\frac{2^N}p 
 \PP\prbig{ \abs{X_0}>
 2^{N/q\pr{ p,\gamma} }  }.
\end{align*}

\section{Proofs}\label{sec:proofs}

\subsection{Proof of Theorem~\ref{thm:Besov_alpha_p_leq_1}}
Similarly as in the proof of \cite[Theorem~4]{MR3736194}, using the continuous
embeddings
$B_{2,\gamma}^o([0,1])\hookrightarrow B_{p,\gamma}^o([0,1])$ if $1\leq p\leq 2$
and
$ \gamma\in \pr{0,1/2}$, it suffices to prove the following cases:
\begin{itemize}
 \item $p=2$ and
$ \gamma\in\pr{0,1/2}$;
\item $p>2$ and $ \gamma\in (0,1/p]$.
\end{itemize}
The convergence of the finite-dimensional distributions of $\pr{W_{n,\mathbb{X}}\pr{\cdot}}_{n\geq 1}$ is already assumed. It
thus suffices to show its tightness, which
translates as
\begin{equation}\label{eq:tightness_p_alpha_leq_1}
 \forall \varepsilon>0,\quad \lim_{\delta\to 0}
 \limsup_{n\to\infty}\PP\prbig{\frac 1{\delta^{p\gamma}}
 \sup_{h:\abs{h}\leq \delta} \int_{I_h}
\abs{W_{n,\mathbb{X}}(t+h)-W_{n,\mathbb{X}}(t)}^p\mathrm{d}t >\varepsilon
  }=0,
\end{equation}
where $I_h=[0,1]\cap [-h, 1,1-h]$. 

We first consider the
case, where $p=2$ and $\gamma\in\pr{0,1/2}$. We bound
the supremum involved in \eqref{eq:tightness_p_alpha_leq_1} by a function of
the squares of increments of the partial sums of the $X_i$. Bounding
them further by maxima of increments and using boundedness in $\el_{2}$ of the
maxima of squares of the normalized partial sums gives
\eqref{eq:tightness_p_alpha_leq_1}.

In the second case, that is, $p>2$  and $\gamma\in (0,1/p]$, the previous argument does not work out
directly. Indeed, we cannot use Markov's inequality, since we do not assume that
the random variables $X_i$ have a finite moment of order $p$. To
overcome this problem, we
will split the probability in \eqref{eq:tightness_p_alpha_leq_1} according to
the cases where $ \displaystyle\sup_{h:\abs{h}\leq \delta} \sup_{t\in
I_h}\abs{W_{n,\mathbb{X}}(t+h)-W_{n,\mathbb{X}}(t)}^{p-2} $ is bigger or not than a fixed
value $\eta$.
The part where this random variable is bigger than $\eta$ is handled by
tightness in $C\pr{[0,1]}$.
On the set where the
previous random variable is smaller,  the $\mathbb L_{p}$-norm involved in
\eqref{eq:tightness_p_alpha_leq_1} reduces to the $\mathbb L_{2}$-norm, for which
the estimates derived in the case $p=2$ apply.

Define the intervals $J_{n,k}:=[(k-1)/n,k/n)$, $k\in\ens{1,\dots,n}$. Let
$\delta,h$ be fixed and such that $\abs{h}\leq\delta$. Let us start by the
case $p=2$ and
 $\gamma\in \pr{0,1/2}$. For $n\geq 1/\delta$,
we bound the random variable
\begin{equation*}
 M\pr{n,\delta}:=\frac 1{\delta^{2\gamma}}
 \sup_{h:\abs{h}\leq \delta} \int_{I_h}
\pr{W_{n,\mathbb{X}}(t+h)-W_{n,\mathbb{X}}(t)}^2\mathrm{d}t,
\end{equation*}
by looking at which interval belong $t$ and $t+h$, so that
\begin{align*}
 \abs{W_{n,\mathbb{X}}(t+h)-W_{n,\mathbb{X}}(t)}^2&=
 \sum_{k,\ell=1}^n \mathbf{1}_{t\in
J_{n,k}} \mathbf{1}_{t+h\in
J_{n,\ell}}\pr{W_{n,\mathbb{X}}(t+h)-W_{n,\mathbb{X}}(t)}^2
\\
 &\leq 3\sum_{k,\ell=1}^n \mathbf{1}_{t\in
J_{n,k}}\mathbf{1}_{t+h\in J_{n,\ell}}
\prbig{W_{n,\mathbb{X}}(t+h)-W_{n,\mathbb{X}}\prbig{\frac{\ell-1}n}}^2\\
&+3\sum_{k,\ell=1}^n \mathbf{1}_{t\in
J_{n,k}}\mathbf{1}_{t+h\in J_{n,\ell}}
\prbig{ W_{n,\mathbb{X}}\prbig{\frac{\ell-1}n}-
W_{n,\mathbb{X}}\prbig{\frac{k-1}n}}^2\\
&+3\sum_{k,\ell=1}^n \mathbf{1}_{t\in
J_{n,k}}\mathbf{1}_{t+h\in J_{n,\ell}}
\prbig{ W_{n,\mathbb{X}}\prbig{\frac{k-1}n}-W_{n,\mathbb{X}}\pr{t}}^2\\
&\leq \frac{2}n\sum_{k=1}^n \mathbf{1}_{t\in J_{n,k}}X_k^2 
+\frac{6}n\sum_{1\leq k,\ell\leq n}\mathbf{1}_{t\in 
J_{n,k}}\mathbf{1}_{t+h\in J_{n,\ell}}
\prbig{\sum_{i=\min\ens{k,\ell}}^{\max\ens{k,\ell}-1}X_i  }^2,
\end{align*}
where we also used the fact that the intervals $J_{n,k}$ are pairwise disjoint in
order to simplify the first and third terms. Taking the integral, letting $J_{n,\ell}-h:=
\ens{u-h\in\R, u\in J_{n,\ell}-h} $, noticing
that the Lebesgue measure of $J_{n_k} \cap \pr{J_{n,\ell}-h}$
is $0$ if $\abs{\ell-k}\geq \ent{n\delta}+2$ and otherwise smaller than
$n^{-1}$ and then using the inequality
\begin{equation*}
 \Big(\sum_{i=\min\ens{k,\ell}}^{\max\ens{k,\ell}-1}X_i\Big)^2
 \leq 4\max_{1\leq j\leq
 2\ent{n\delta}}\Big(\sum_{i=k-\ent{n\delta}+1}^{k-\ent{n\delta}+j}
X_i\Big)^2,
\end{equation*}
 we derive that
\begin{align*}
 M\pr{n,\delta}&\leq \frac 2{n^2\delta^{2\gamma}}
 \sum_{k=1}^n X_k^2 
 +\frac{24}{n\delta^{2\gamma}}\sum_{1\leq k,\ell\leq n}
 \lambda\pr{J_{n_k} \cap \pr{J_{n,\ell}-h}}
\max_{1\leq j\leq
 2\ent{n\delta}}\prbig{\sum_{i=k-\ent{n\delta}+1}^{k-\ent{n\delta}+j}
X_i}^2\\
&\leq \frac 2{n^2\delta^{2\gamma}}
 \sum_{k=1}^n X_k^2 
 +\frac{24}{n^2\delta^{2\gamma}}\sum_{k=1}^n
\max_{1\leq j\leq
 2\ent{n\delta}}\prbig{\sum_{i=k-\ent{n\delta}+1}^{k-\ent{n\delta}+j}
X_i}^2.
\end{align*}
As a consequence, we set
\begin{equation*}
 A:=\sup_{n\geq 1}\frac 1n\mathbb{E}\Big[\max_{1\leq k\leq n}
    S_{k,X_0}^2\Big],
\end{equation*}
which is finite by \eqref{eq:uniform_integrability_for_p_alpha_leq_1}, we infer
that
\begin{align}
 \PP\pr{M\pr{n,\delta}>\eps}&\leq \frac 1{\eps}
 \E{M\pr{n,\delta}}\nonumber\\
 &\leq \frac 2{n\delta^{2\gamma}\eps }\mathbb{E}\big[X_0^2\big]
+\frac{24}{n\delta^{2\gamma}\eps}\sum_{k=1}^n
\mathbb{E}\Big[
\max_{1\leq j\leq
 2\ent{n\delta}}\prbig{\sum_{i=k-\ent{n\delta}+1}^{k-\ent{n\delta}+j}
X_i}^2\Big]\nonumber \\
&\leq  \frac 2{n\delta^{2\gamma}\eps }\mathbb{E}\big[X_0^2\big]
+\frac{48}{n^2\delta^{2\gamma}\eps}n \ent{2n\delta}
A\nonumber
\\
&\leq \frac 2{n\delta^{2\gamma}\eps
}\mathbb{E}\big[X_0^2\big]+\frac{96}\eps A\delta^{1-2\gamma},\label{eq:bound_tail_Mndelta}
\end{align}
from which \eqref{eq:tightness_p_alpha_leq_1} follows.

Let us now treat the case $p>2$ and $\gamma\in (0,1/p]$. For fixed $n\geq 1$ and positive
$\delta$, $\eta$, we define the event
\begin{equation*}
 A\pr{n,\delta,\eta}:=\Big\{\sup_{h:\abs{h}\leq \delta} \sup_{t\in I_h} 
\abs{W_{n,\mathbb{X}}(t+h)-W_{n,\mathbb{X}}(t)}^{p-2}> \eta\Big\}.
\end{equation*}
Condition \eqref{eq:uniform_integrability_for_p_alpha_leq_1} implies tightness
of $\pr{W_{n,\mathbb{X}}\pr{\cdot}}_{n\geq 1}$ in
$C\pr{[0,1]}$. Therefore, for each positive $\eta$,
\begin{equation*}
    \lim_{\delta\to 0}\limsup_{n\to\infty}\PP\pr{A\pr{n,\delta,\eta}}=0.
\end{equation*}
Moreover, the inclusion
\begin{multline*}
\Big\{\frac 1{\delta^{p\gamma}}
 \sup_{h:\abs{h}\leq \delta} \int_{I_h}
  \abs{W_{n,\mathbb{X}}(t+h)-W_{n,\mathbb{X}}(t)} ^p\mathrm{d}t
>\varepsilon
  \Big\}\cap A\pr{n,\delta,\eta}^c\\
\subset  \Big\{\frac 1{\delta^{p\gamma}}
 \sup_{h:\abs{h}\leq \delta} \int_{I_h}
\pr{  W_{n,\mathbb{X}}(t+h)-W_{n,\mathbb{X}}(t)} ^2\mathrm{d}t
>\frac{\varepsilon}{\eta}
  \Big\}
=\Big\{\delta^{\pr{2-p}\gamma} M\pr{n,\delta} > \frac{\varepsilon}{\eta} \Big\}
\end{multline*}
takes place. Applying \eqref{eq:bound_tail_Mndelta}
with $\eps$ replaced by $\eps \eta^{-1}\delta^{\pr{p-2}\gamma}$ gives the upper bound 
\begin{equation*}
 \PP\prbig{\frac 1{\delta^{p\gamma}}
 \sup_{h:\abs{h}\leq \delta} \int_{I_h}
  \abs{W_{n,\mathbb{X}}(t+h)-W_{n,\mathbb{X}}(t) }^p\mathrm{d}t
>\varepsilon}  
\leq  \PP\pr{A\pr{n,\delta,\eta}}+\frac{2}{n\delta^{p\gamma}}\mathbb{E}\big[X_0^2\big]+
\frac{96}{\eps}A\delta^{1-p\gamma}\eta.
\end{equation*}
Hence, taking $\displaystyle\limsup_{\delta\to 0}\limsup_{n\to\infty}$ ensures that 
for each positive $\eta$,
\begin{equation*}
 \limsup_{\delta\to 0}\limsup_{n\to\infty}\PP\prbig{\frac 1{\delta^{p\gamma}}
 \sup_{h:\abs{h}\leq \delta} \int_{I_h}
 \abs{W_{n,\mathbb{X}}(t+h)-W_{n,\mathbb{X}}(t)}^p\mathrm{d}t
>\varepsilon}
\leq
\frac{96}{\eps}A \eta\mathbf{1}_{p\gamma=1 },
\end{equation*}
which proves \eqref{eq:tightness_p_alpha_leq_1} and concludes the proof of Theorem~\ref{thm:Besov_alpha_p_leq_1}.

\subsection{Proof of Theorem~\ref{thm:Hannan}}

The convergence of finite-dimensional distributions follows from \cite[Corollary 2]{MR2359065}, which establishes the WIP in $C([0,1])$ under the Hannan condition \eqref{eq:Hannan}.

In order to prove tightness of $\pr{W_{n,\mathbb{X}}\pr{\cdot}}_{n\geq 1}$ in $B_{p,\gamma}^o\pr{[0,1]}$, we use  \cite[Theorem~8]{MR3736194}, which reduces to check 
\begin{equation}\label{eq2.1}
 \lim_{R\to\infty}\sup_{n\geq 1}\PP\prbig{\norm{W_{n,\mathbb{X}}
  \pr{\cdot}}_{\el_p\pr{[0,1]}}  >R}=0 \text{ and}
\end{equation}
\begin{equation}\label{eq2.2}
\forall \eps>0,\ \  \lim_{\delta\to 0}\sup_{n\geq 1}
 \PP\prbig{\frac 1{\delta^\gamma}\omega_p\pr{W_{n,\mathbb{X}}
 \pr{\cdot},\delta}>\eps }=0.
\end{equation}
A possibility would be to apply Proposition~\ref{prop:tightness_criterion} in order to check tightness. To do so, one would need a deviation inequality, which is, to the best of our knowledge, not available in the literature. Moreover, the decomposition of $X_i$ as a sum of projectors would lead to a sub-optimal result since the required condition would read as $\sum_{i\in\Z}\norm{\E{X_{0}\mid\Fca_{i}}-
\E{X_{0}\mid\Fca_{i-1}}  }_{q\pr{p,\gamma}}^{s}<\infty$ for some $s \in (0,1)$.

Let us show \eqref{eq2.1}. It follows from \cite[Corollary~2]{MR2359065}
 that $\pr{W_{n,\mathbb{X}}\pr{\cdot}}_{n\geq 1}$ is tight in $C([0,1])$. In particular, $\pr{\|W_{n,\mathbb{X}}\pr{\cdot}\|_{\infty}}_{n\geq 1}$ is tight.
Since $\|x\|_{\el_p\pr{[0,1]}}\leq \|x\|_{\infty}$ for each $x\in C([0,1])$ and $p\geq 1$, we get the desired relation \eqref{eq2.1}.

It remains to establish \eqref{eq2.2}. To this aim, we first decompose $X_{0}$ as $\sum_{i \in \mathbb{Z}} X_{0}^{(i)}$, where for each $i\in\Z$, 
$$
X_{0}^{(i)}:=\mathbb{E}[X_0 \mid \mathcal{F}_i]-\mathbb{E}[X_0 \mid \mathcal{F}_{i-1}], 
\quad \mathbb{X}^{(i)}:=\big(X_{0}^{(i)}\circ \mathbb{T}^k\big)_{k\in\Z} .
$$ 
Then, for each fixed non-negative integer $L$, we have the upper bound
$$
\PP\pr{\frac 1{\delta^\gamma}\omega_p\pr{W_{n,\mathbb{X}}
 \pr{\cdot},\delta}>\eps } \leq P_{1}(n,\delta,L)+P_{2}(n,\delta,L),$$
where
\begin{align*}
   P_{1}\pr{n,\delta,L}&:=
 \PP\prbig{\frac 1{\delta^\gamma}\omega_p\prbig{  \sum_{|i|\leq L }  W_{n,\mathbb{X}^{(i)}}
 \pr{\cdot},\delta}>\frac{\eps}{2} },\\
P_{2}\pr{n,\delta,L}&:=\PP\prbig{\frac 1{\delta^\gamma}\omega_p\prbig{ \sum_{|i|> L }W_{n,\mathbb{X}^{(i)}}
 \pr{\cdot},\delta}>\frac{\eps}{2} }. 
\end{align*}
To prove that \eqref{eq2.2} holds true, it suffices to prove that, for any $\varepsilon>0$,
\begin{equation}\label{eq:Hannan_contribution_finite_sum}
\forall L\geq 1,\quad \lim_{\delta\to 0}\sup_{n\geq 1}P_{1}(n,\delta,L)=0,
\end{equation}
\begin{equation}\label{eq:Hannan_contribution_finite_tail}
\lim_{L\to\infty}\lim_{\delta\to 0}\sup_{n\geq 1}P_{2}(n,\delta,L)=0.
\end{equation}
Consider $P_{1}(n,\delta,L)$ for fixed $L\geq 1$. Applying the triangle inequality for the modulus of regularity $\omega_{p}(\cdot,\delta)$ then using the inequality 
\begin{equation*}
    \PP\prbig{\sum_{i\in I} Z_i>t}\leq 
 \sum_{i\in I}\PP\prbig{Z_i>\frac{t}{\operatorname{Card}\pr{I}}},
\end{equation*}
valid for each non-empty finite set $I$, $t>0$ and every  non-negative random variables $\pr{Z_i}_{i\in I}$, one has
\begin{align*}
    P_{1}(n,\delta,L)&\leq \PP\prbig{\frac 1{\delta^\gamma}\sum_{|i|\leq L}\omega_p\pr{ W_{n,\mathbb{X}^{(i)}}
 \pr{\cdot},\delta}>\frac{\eps}{2} }\\
 &\leq (2L+1)\sum_{|i|\leq L} \PP\prbig{\frac 1{\delta^\gamma}\omega_p\big( W_{n,\mathbb{X}^{(i)}}
 \pr{\cdot},\delta \big)>\frac{\eps}{2(2L+1)}}. 
\end{align*}
Note that for each fixed $i$, the sequence $\mathbb{X}^{(i)}$ is a martingale difference sequence, where $X_0^{\pr{i}}\in\el_{q\pr{p,\gamma}}$. An application of \cite[Theorem~1]{MR3736194} 
yields tightness of 
$\big(W_{n,\mathbb{X}^{(i)}}\pr{\cdot}\big)_{n\geq 1}$ in $B_{p,\gamma}^o\pr{[0,1]}$. Since $L$ is a fixed, \eqref{eq:Hannan_contribution_finite_sum} follows.

We now consider $P_{2}(n,\delta,L)$. Bounding the modulus of regularity by the Besov norm, then  using the triangle and Markov's inequalities provide the bounds
\begin{equation}\label{eq:bound_P2nL_Hannan}
    P_{2}(n,\delta,L) \leq \PP\prbig{  \sum_{|i|>L} \|W_{n,\mathbb{X}^{(i)}}\pr{\cdot}  \|_{B_{p,\gamma}^{o}\pr{[0,1]}}> \frac{\varepsilon}{2}} \leq \frac{2}{\varepsilon}\sum_{|i|>L} \big\lVert \|W_{n,\mathbb{X}^{(i)}}\pr{\cdot} \|_{B_{p,\gamma}^{o}\pr{[0,1]}}\big\rVert_{1}.
\end{equation}
We now prove that there exists a positive constant $C_{p,\gamma}$ such that for any stationary martingale difference sequence $\mathbb{D}= (D_{i})_{i \in \mathbb{Z}}$, the following inequality holds:
\begin{equation}\label{eq:ineq_Besov_norm_partial_sum_process}
\sup_{n\geq 1}  \left\| \|W_{n,\mathbb{D} }(\cdot) \|_{B_{p,\gamma}^o\pr{[0,1]}} \right\|_1 \leq C_{p,\gamma} \|D_0\|_{q(p,\gamma)},
\end{equation}
where $W_{n,\mathbb{D} }\pr{\cdot}$ is defined in \eqref{eq:def_partial_sum_process} with $X_0$ replaced by $D_0$.\\
Taking into account \eqref{sequential norm}, for any $x>0$, let us first consider the probability
$$
\PP\pr{  \left  \|W_{n,\mathbb{D}}(\cdot) \right \|^{\operatorname{seq}}_{B_{p,\gamma}^o\pr{[0,1]}}>x }.
$$
Observe that
$$
\PP\pr{ \left  \|W_{n,\mathbb{D}}(\cdot) \right \|^{\operatorname{seq}}_{B_{p,\gamma}^o\pr{[0,1]}}>x }
\leq I_{1}(n,x)+I_{2}(n,x),$$
where
$$ I_{1}(n,x):=\PP\prbig{ \sup_{j\leq \lfloor\log_2 n\rfloor} 2^{j(\gamma-1/p)}\prbig{\sum_{r\in \mathbf{D}_{j}}\left|\lambda_{r}(W_{n,\mathbb{D}}(\cdot))\right|^{p} }^{1/p}  >x },
$$
$$ I_{2}(n,x):=
\PP\prbig{ \sup_{j> \lfloor\log_2 n\rfloor} 2^{j(\gamma-1/p)}\prbig{\sum_{r\in \mathbf{D}_{j}}\left|\lambda_{r}(W_{n,\mathbb{D}}(\cdot))\right|^{p} }^{1/p}  >x },
$$
and $\log_2$ denotes the binary logarithm.\\
We are reduced to show the existence of a positive constant $C_{p,\gamma}$ such that 
for any strictly stationary martingale difference sequence $\pr{D_i}_{i\in\Z}$, 
\begin{equation}\label{ineq1}
  \int_0^\infty  I_{1}(n,x)\mathrm{d}x \leq C_{p,\gamma} \|D_0\|_{q(p,\gamma)}
\end{equation}
and 
\begin{equation}\label{ineq2}
    \int_0^\infty  I_{2}(n,x)\mathrm{d}x \leq C_{p,\gamma} \|D_0\|_{q(p,\gamma)}.
\end{equation}
We prove \eqref{ineq1}. Using a union bound, one has 
\begin{equation}\label{eq:union_bound_I1}
    I_{1}(n,x)\leq \sum_{j=1}^{\ent{\log_2n}}\PP\prbig{  2^{j(\gamma-1/p)}\prbig{\sum_{r\in \mathbf{D}_{j}}\left|\lambda_{r}(W_{n,\mathbb{D}}(\cdot))\right|^{p}}^{1/p}  >x }.
\end{equation}
Let $N$ be such that $2^N\leq n<2^{N+1}$ and $j$ be a fixed positive integer such that 
$j\leq \ent{\log_2 n}$.
We now apply the next elementary inequality (obtained via a truncation argument), valid for arbitrary non-negative random variables $Z_1,\ldots,Z_m$, $m,p\geq 1$ and $t>0$, to the right-hand side of \eqref{eq:union_bound_I1}:
\begin{equation}\label{eq:lp-tail}
 \mathbb P\prbig{
 \prbig{\sum_{i=1}^m Z_i^p}^{1/p}>t}
 \leq C_p\sum_{i=1}^m
 \int_{0}^{1}u^{p-1}\mathbb P(Z_i>tu)\,\mathrm{d}u.  
\end{equation}
Hence, we obtain 
\begin{equation*}
    I_{1}(n,x)\leq C_p\sum_{j=1}^{\ent{\log_2n}}
    \sum_{r\in\mathbf{D}_j}\int_0^1 u^{p-1}\PP\prbig{\abs{\lambda_r\pr{W_{n,\mathbb{D}}
    \pr{\cdot} }}>xu
     2^{j\pr{1/p-\gamma}} }
    \mathrm{d}u.
\end{equation*}
For $r\in \mathbf{D}_j$, $j\leq \ent{\log_2 n}$, the increment $\lambda_r\pr{S_{n,\mathbb{D}}}$ is a convex
combination of discrete increments of the partial sums, each involving at
most $C2^{N-j}$ consecutive variables, where $C$ is a numerical constant. Hence, by stationarity of $\mathbb{D} $ and from Proposition \eqref{devationinequality}, we infer that for $q:=\pr{q\pr{p,\gamma}+p}/2$,
$$
\mathbb P(\abs{\lambda_r\pr{W_{n,\mathbb{D}}\pr{\cdot} }}>t)
 \leq C_{p,\gamma}\int_0^\infty
 \mathbb P\left(
 |D_0|>c\,t\,2^{-(N-j)/2}v
 \right)
 \min\{v^{q-1},v\}\,\mathrm{d}v .
$$
Combining\eqref{eq:lp-tail} together with $\operatorname{Card}\pr{\mathbf{D}_j}=2^{j-1}$, it follows that 
\begin{equation*} 
I_1\pr{n,x}
\leq C_{p,\gamma} \sum_{j=1}^{\ent{\log_2 n}}2^j\int_0^1\int_0^\infty\min\{v^{q-1},v\}\,u^{p-1}
 \mathbb P\left(
 |D_0|>xuv\,2^{j/q\pr{p,\gamma}}
 \right)
 \mathrm{d}v\,\mathrm{d}u ,
\end{equation*}
where we have used the fact that
$$
\frac{\sqrt n\,2^{-j(\gamma-1/p)}}{2^{(N-j)/2}}
\leq  
2^{j(1/2-(\gamma-1/p)+1/2}
= 2^{j/q\pr{p,\gamma}+1/2},
$$
and the constant $C_{p,\gamma}$ may change from line to line and depends only on $p$ and $\gamma$. \\
In view of the following inequality
$$
\sum_{j=0}^{\infty} 2^j
\mathbb P(Z>a2^{j/q\pr{p,\gamma}})
\leq
C_{p,\gamma} a^{-q\pr{p,\gamma}}\mathbb{E}\big[Z^{q\pr{p,\gamma}}\big],
$$
which holds for any non-negative random variable $Z$ and $a>0$, we claim that
$$
 I_1\pr{n,x}
\leq
C x^{-q\pr{p,\gamma}}\E{|D_0|^{q\pr{p,\gamma}}}
\int_0^1u^{p-q\pr{p,\gamma}-1}\,\mathrm{d}u
\int_0^\infty
v^{-q\pr{p,\gamma}}\min\{v^{q-1},v\}\,\mathrm{d}v.
$$
Since $2<q(p,\gamma)<q<p$, the integrals over $u$ and $v$ are finite and 
$I_1\pr{n,x}$ is obviously bounded by $1$, one arrives at  
\begin{equation*}
    I_1\pr{n,x}
\leq
\min\big\{
1,\,
C_{p,\alpha}x^{-q(p,\gamma)}\|D_{0}\|_{q(p,\gamma)}^{q(p,\gamma)}
\big\}.
\end{equation*}
Integrating this estimate over $x>0$ gives us the desired bound \eqref{ineq1}. 

It remains to prove \eqref{ineq2}.
We follow the argument of \cite{MR3736194} (see the top of page 449) for the dyadic levels such that if one belongs to an interval of the form $[k/n,\pr{k+1}/n)$, the next one belongs to it or to $[(k+1)/n,\pr{k+2}/n)$. 
For $j>\lfloor\log_2 n\rfloor$, using the polygonal structure of
$W_{n,\mathbb{D}}\pr{\cdot}$ and the corresponding cardinality estimates, we obtain 
the pointwise inequality
$$
\sup_{j>\lfloor\log_{2} n\rfloor}
2^{j(\gamma-1/p)}
\prbig{
\sum_{r\in \mathbf{D}_j}
\abs{\lambda_r\pr{W_{n,\mathbb{D}}\pr{\cdot}}}^p
}^{1/p}
\leq
C_{p,\gamma} n^{-1/q\pr{p,\gamma}}
\prbig{
\sum_{k=1}^n |D_k|^p
}^{1/p},
$$
where we have used that  $
-1/2+\gamma-1/p=-1/q(p,\gamma)$.

Consequently,
$$
I_2\pr{n,x}
\leq
\mathbb{P}\prbig{
\prbig{\sum_{k=1}^n|D_k|^p}^{1/p}
>\frac{x n^{1/q(p,\gamma)}}{C_{p,\gamma}}
}.
$$
It follows from \eqref{eq:lp-tail} and the stationarity of $(D_k)_{k\in\Z}$  that
$$
I_2\pr{n,x}
\leq
C_{p,\gamma} n\int_{0}^{1}u^{p-1}
\mathbb{P}\left(
|D_0|> x n^{1/q(p,\gamma)}u
\right)\,\mathrm{d}u.
$$
Therefore, by Markov's inequality at order $q(p,\gamma)$ and in view of $p>q$, 
$$
I_2(n,x)
\leq
C_{p,\gamma} x^{-q(p,\gamma)}\norm{D_0}_{q(p,\gamma)}^{q(p,\gamma)}
\int_0^1 u^{p-q-1}\,\mathrm{d}u.
$$
Moreover, we also have $I_2\pr{n,x}\leq 1$ and hence 
\begin{equation*}
    I_2\pr{n,x}\leq \min\big\{1,C_{p,\gamma}x^{-q(p,\gamma)}\norm{D_0}_{q(p,\gamma)}^{q(p,\gamma)}\big\}.
\end{equation*}
Integrating this estimate over $x>0$ concludes the proof of \eqref{eq:ineq_Besov_norm_partial_sum_process}. 

Overall, combining \eqref{eq:bound_P2nL_Hannan} with \eqref{eq:ineq_Besov_norm_partial_sum_process}, 
we find that for some positive constant $C_{p,\gamma}$, 
\begin{equation*}
\sup_{\delta\in\pr{0,1}}  \sup_{n\geq 1} P_{2}(n,\delta,L) \leq 
\frac{C_{p,\gamma}}{\varepsilon}\sum_{|i|>L} \lVert X_{0}^{(i)}\rVert_{q\pr{p,\gamma}}
=\frac{C_{p,\gamma}}{\varepsilon}\sum_{|i|>L} \norm{\E{X_0\mid\Fca_{i}}
-\E{X_0\mid\Fca_{i-1}}}_{q\pr{p,\gamma}}.
\end{equation*} 
As $L \to \infty$, the latter sum vanishes by \eqref{eq:Hannan}. Thus, \eqref{eq:Hannan_contribution_finite_tail} holds, which ends the proof of Theorem~\ref{thm:Hannan}.
\subsection{Proof of Corollary~\ref{cor:Bernoulli}}
For $i\in\Z$, let $\Fca_i:=\sigma\pr{\eps_j, j\leq i}$.
It suffices to see that condition \eqref{eq:Hannan} implies
\eqref{eq:cond_Bernoulli}, which is a consequence of the bound 
$$\delta_{\mathbb{X},i,q\pr{p,\gamma}}\leq  \norm{\E{X_i\mid\Fca_{0}}-
\E{X_i\mid\Fca_{-1}}  }_{q\pr{p,\gamma}}= \norm{\E{X_0\mid\Fca_{-i}}-
\E{X_0\mid\Fca_{-i-1}}  }_{q\pr{p,\gamma}},$$
where the inequality comes form \cite[Theorem~1 (i)]{MR2172215} and the equality is a direct consequence of strict stationarity of $\mathbb{X}$.
\subsection{Proof of Theorem~\ref{thm:alpha_dep}}
The convergence of the finite-dimensional distributions of $\pr{W_{n,\mathbb{X}}\pr{\cdot}}_{n\geq 1}$ is addressed  by \cite[Theorem~3.3]{MR3583992}.
It remains to verify its tightness. To this aim, we have to check
\eqref{eq:conv_for_tightness}.
For fixed $N$, $j$ and $u\in [0,1]$, we apply Proposition~\ref{prop:deviation_alpha_dep} with 
$\ell:=N-j$, $x:=u2^{N/2}2^{j\pr{1/p-\gamma}}$, $r \in (2q(p,\gamma)-2, p)$ and $\beta \in (r-2,r)$ (which is possible by the assumption that $\gamma <1/2+1/p-2/\pr{p+2}$).

We derive the bound
$$
\sum_{j=J}^{N} 2^{j} \int_{0}^{1} u^{p-1} \mathbb{P}\left(\max_{1 \leq k \leq 2^{N-j}} 
\abs{S_{k,X_0}} > 2^{N/2} 2^{j(1/p-\gamma)} u \right) \mathrm{d}u \leq C_{p,\gamma,r} \sum_{i=1}^4  I_{i}^{\alpha}\pr{J,N},
$$
where 
\begin{align*}
 I_1^{\alpha}\pr{J,N}&:=   \sum_{j=J}^N2^j
 \int_0^1  u^{p-1}\frac{2^{(N-j)r/2}}{(2^{N/2}2^{j(1/p - \gamma)} u)^r} 
 \prbig{\sum_{k=0}^{2^{N-j}-1}\int_0^{\alpha_{1,Y}(k)}Q^{2}(v)\mathrm{d}v}^{r/2}\mathrm{d}u, \\
 I_2^{\alpha}\pr{J,N}&:= \sum_{j=J}^N2^j \int_0^1  u^{p-1}\frac{2^{N-j}}{2^{N/2}2^{j(1/p - \gamma)} u}\int_0^{L (2^{N/2}2^{j(1/p-\gamma)}u)}Q(v)\mathrm{d}v\mathrm{d}u,\\
 I_3^{\alpha}\pr{J,N}&:= \sum_{j=J}^N2^j \int_0^1  u^{p-1}\frac{2^{N-j}}{(2^{N/2}2^{j(1/p -\gamma)} u)^{1+\beta/2}}\int_0^{L (2^{N/2}2^{j(1/p - \gamma)} u)}R_{\pr{\alpha}}^{\beta/2}(v)Q(v)\mathrm{d}v\mathrm{d}u,\\
  I_4^{\alpha}\pr{J,N}&:= \sum_{j=J}^N2^j \int_0^1  u^{p-1}\frac{2^{N-j}}{(2^{N/2}2^{j(1/p-\gamma)} u)^{1+p/2}}\int_{L (2^{j(1/p-\gamma)} u)}^{1}R_{\pr{\alpha}} ^{r/2}(v)Q(v)\mathrm{d}v\mathrm{d}u.
\end{align*}
We have to establish that for each $i\in\ens{1,2,3,4}$, 
\begin{equation}\label{eq:conv_I_i_alpha}
    \lim_{J \to \infty} \limsup_{N \to \infty} 
    I_{i}^\alpha\pr{J,N}=0.
\end{equation}
 Let us bound $I_1^{\alpha}\pr{J,N}$. Bounding the sum over $k$ and the integral over $v$ by $C:=\int_0^1 \alpha_{1,Y}^{-1}(v) Q^2(v) \, \mathrm{d}v < \infty$, we infer that
$$
I_1^{\alpha}\pr{J,N} \leq C^{r/2} \int_{0}^1u^{p-r-1}\mathrm{d}u\sum_{j=J}^{\infty}
2^{j\pr{1-r/q\pr{p,\gamma}}},
$$
hence, since $r>q\pr{p,\gamma}>2$, \eqref{eq:conv_I_i_alpha} holds for $i=1$.
 
We now consider $I_{2}^{\alpha}\pr{J,N}$. We define the function $\eps_1^{\alpha}\colon \R_+\to\R$ by 
\begin{equation*}
    \eps_1^{\alpha}\pr{t}:=\sup_{x\geq t}x^{q\pr{p,\gamma}-1}\int_0^1Q\pr{v}\mathbf{1}_{R_{\pr{\alpha}}\pr{v}\geq x}\mathrm{d}v.
\end{equation*}
By assumption \eqref{eq:cond_alpha_dep}, $\displaystyle\lim_{t\to\infty}\eps_1^{\alpha}\pr{t}=0$. 
Using properties of inverse functions, one can rewrite the integral over $v$ as $\int_0^1 Q\pr{v}
\mathbf{1}_{R_{\pr{\alpha}}\pr{v}> 2^{j(1/p-\gamma)}u }\mathrm{d}v$ and then bounding it using the function $\eps_1^{\alpha}\pr{\cdot}$ yields
\begin{align*}
    I_2^\alpha\pr{J,N}&\leq \sum_{j=J}^N2^{N/2}2^{j\pr{\gamma-1/p}}\int_0^1u^{p-2}
    \prbig{ 2^{N/2}2^{j(1/p-\gamma)}u}^{1-q\pr{p,\gamma}}\eps_1^{\alpha}\prbig{ 2^{N/2}2^{j(1/p-\gamma)}u}\mathrm{d}u\\
    &\leq 2^{N\pr{1-q\pr{p,\gamma}/2 }}
    \sum_{j=1}^N2^{j\pr{1/p-\gamma}q\pr{p,\gamma}-1}
    \int_0^1u^{p-q\pr{p,\gamma}}\eps_1^{\alpha}\prbig{2^{N/q\pr{p,\gamma}}u }\mathrm{d}u\\
    &\leq C_{p,\gamma}\int_0^1u^{p-q\pr{p,\gamma}}\eps_1^{\alpha}\prbig{2^{N/q\pr{p,\gamma}}u }\mathrm{d}u,
\end{align*}
and it follows by the dominated convergence theorem 
that \eqref{eq:conv_I_i_alpha} is also valid for $i=2$. 

Consider $I_{3}^{\alpha}(J,N)$. 
Define the function $\eps_2^{\alpha}\colon \R_+\to\R$ by 
\begin{equation*}
    \eps_2^{\alpha}\pr{t}:=\sup_{x\geq t}x^{q\pr{p,\gamma}-1-\frac{\beta}2}\int_0^1R_{\pr{\alpha}}^{\beta/2}\pr{u}\mathbf{1}_{R\pr{u}>x}\mathrm{d}u.
\end{equation*}
By \cite[Lemma~6.2]{MR3583992}, one has that $\displaystyle\lim_{t\to\infty}\eps_2^{\alpha}\pr{t}=0$.
Proceeding as before, that is, estimating the integral over $v$ using the function $\eps_2^{\alpha}\pr{\cdot}$, then majoring the argument of $\eps_2^{\alpha}\pr{\cdot}$ uniformly over $j$ and summing over $j$ gives
$$ I_{3}^{\alpha}\pr{J,N} \leq C_{p,\gamma}  \int_0^1u^{p-q\pr{p,\gamma}}\eps_2^{\alpha}\prbig{2^{N/q\pr{p,\gamma}}u }\mathrm{d}u,$$
hence by the dominated convergence theorem, \eqref{eq:conv_I_i_alpha} holds true for $i=3$.

It remains to consider $I_{4}^{\alpha}\pr{J,N}$.  
Define the function $\eps_3^{\alpha}\colon \R_+\to\R$ by 
\begin{equation*}
    \eps_3^{\alpha}\pr{t}:=\sup_{x\geq t}x^{-\frac r2+q\pr{p,\gamma}-1}\int_0^1R_{\pr{\alpha}}^{r/2}\pr{u}\mathbf{1}_{R\pr{u}>x}\mathrm{d}u.
\end{equation*}
By \cite[Lemma~6.2]{MR3583992}, one has that $\displaystyle\lim_{t\to\infty}\eps_3^{\alpha}\pr{t}=0$.
The same arguments as those used for bounding the terms $I_2^{\alpha}\pr{J,N}$ and $I_3^{\alpha}\pr{J,N}$
yield the bound 
$$ I_{4}^{\alpha}\pr{J,N} \leq C_{p,\gamma} 
\int_0^1u^{p-q\pr{p,\gamma}}\eps_3^{\alpha}\prbig{2^{N/q\pr{p,\gamma}}u }\mathrm{d}u.$$
Again by the dominated convergence theorem, \eqref{eq:conv_I_i_alpha} is also valid for $i=4$. This ends the proof of Theorem~\ref{thm:alpha_dep}.
\subsection{Proof of Theorem~\ref{thm:rho_mixing}}
The convergence of the finite-dimensional distributions of $\left(W_{n,\mathbb{X}}\pr{\cdot}\right)_{n \geq 1}$ has been addressed in \cite[Corollary~15]{MR2206313}. Thus, the proof will be finished if we verify tightness of $\left(W_{n,\mathbb{X}}\pr{\cdot} \right)_{n \geq 1}$ in $B_{p,\gamma}^{o}\pr{[0,1]}$, which will be done by checking \eqref{eq:conv_for_tightness}.

We  apply Proposition~\ref{prop:deviation_rho_mixing}  with $\ell:=N-j$,
$x:=u2^{N/2}2^{j\pr{1/p-\gamma}}$,
$A:=\eta u^{\frac 1{1-q\pr{p,\gamma}}} 2^{N/2}2^{j\pr{1/p-\gamma}}$, $q\in \pr{q\pr{p,\gamma},p}$ and $\eta$ is an arbitrary but fixed positive number.
We first need to make sure that $A$ is indeed in the range of application of Proposition~\ref{prop:deviation_rho_mixing}, that is, there exists a positive integer $J_{\eta}$ such that for each $N\geq J\geq J_{\eta}$,
$$\forall j \text{ such~that } J\leq j\leq N,\quad 2\cdot2^{N-j}\mathbb{E}[|X_{0}|\mathbf{1}_{|X_{0}|\geq A}] \leq x.$$
Denoting
\begin{equation}\label{eq:def_M_weal_l_q_norm}
M:=\ \sup_{t>0}t^{q\pr{p,\gamma}}\PP\pr{\abs{X_0}
>t } ,
\end{equation}
we obtain the bound 
\begin{align*}
 \E{|X_{0}|\mathbf{1}_{|X_{0}|\geq A}}&=
A\PP\pr{\abs{X_0}\geq A}+\int_{A}^{\infty}\mathbb{P}(|X_{0}|\geq t)
\mathrm{d}t
\\
&\leq A^{1-q\pr{p,\gamma}}M+
M\int_{A}^{\infty}\frac{\mathrm{d}t}{t^{q(p,\gamma)}}\\
&=M A^{1-q\pr{p,\gamma}}\frac{q\pr{p,\gamma}}{q\pr{p,\gamma}-1}.
\end{align*}
Therefore,
\begin{align*}
\frac{2}{x}\cdot2^{N-j}\mathbb{E}[|X_{0}|\mathbf{1}_{|X_{0}|\geq A}]&\leq
\frac{M}{x} A^{1-q\pr{p,\gamma}}\frac{q\pr{p,\gamma}}{q\pr{p,\gamma}-1}\\
&=\frac{2Mq\pr{p,\gamma}}{q\pr{p,\gamma}-1}\eta^{1-q\pr{p,\gamma}}
2^{N\pr{1-q\pr{p,\gamma}/2}}2^{-j\pr{1+q\pr{p,\gamma}\pr{1/p-\gamma}  }  } \\
&=\frac{2Mq\pr{p,\gamma}}{q\pr{p,\gamma}-1}\eta^{1-q\pr{p,\gamma}}
2^{N\pr{1-q\pr{p,\gamma}/2}}2^{-j\pr{2-q\pr{p,\gamma}/2}}\\
&\leq 2^{-J}\frac{2Mq\pr{p,\gamma}}{q\pr{p,\gamma}-1}\eta^{1-q\pr{p,\gamma}},
\end{align*}
which is smaller than $1$ for all $J\geq J_\eta$.

The application of inequality \eqref{eq:deviation_rho_mixing} gives
\begin{equation}\label{eq:decomposition_proba_rho_mixing}
\sum_{j=J}^{N}2^{j}\int_{0}^{1}u^{p-1}\mathbb{P}\prbig{\max_{1\leq k\leq 2^{N-j}}\abs{S_{k,X_{0}}}> u2^{N/2}2^{j\pr{1/p -\gamma}} }\mathrm{d}u 
\leq I_{1}^{\rho}\pr{J,N,\eta}+I_{2}^{\rho}\pr{J,N}+I_{3}^\rho\pr{J,N,\eta},
\end{equation}
where
\begin{align*}
 I_{1}^{\rho}\pr{J,N,\eta}&:=\sum_{j=J}^{N}2^{j}\int_{0}^{1}u^{p-1}
 2^{N-j}\mathbb{P}\pr{\abs{X_{0}}\geq \eta u^{1/\pr{1-q\pr{p,\gamma}}}
2^{N/2}2^{j\pr{1/p-\gamma}}}\mathrm{d}u, \\
 I_{2}^{\rho}(J,N)&:=K\sum_{j=J}^{N}2^{j}\int_{0}^{1}u^{p-1}(u2^{j/q(p,\gamma)}2^{(N-j)/2})^{-q}2^{(N-j)q/2}\exp\Big(K\sum_{i=1}^{N-j}\rho_{\mathbb{X}}(2^{i})\Big)
\|X_{0}\|_{2}^{q}\mathrm{d}u, \\
I_{3}^{\rho}\pr{J,N,\eta}&:=K\sum_{j=J}^{N}2^{j}\int_{0}^{1}u^{p-1}(u2^{j/q(p,\gamma)}2^{(N-j)/2})^{-q}2^{N-j}\exp\prbig{\sum_{i=1}^{N-j}\left(\rho_{\mathbb{X}}(2^{i})
\right)^{2/q}}\|X_{0}\mathbf{1}_{|X_{0}|\leq A}\|_{q}^{q}\mathrm{d}u.
\end{align*}
In order to  estimate $I_{1}^{\rho}\pr{J,N,\eta}$, we define the function
$\eps^{\rho}\colon \R_+\to\R$ by 
 \begin{equation*}
 \eps^{\rho}\pr{t}:=\sup_{y\geq t}y^{q\pr{p,\gamma}}\PP\pr{\abs{X_0}>y}.
 \end{equation*}
By assumption \eqref{eq:tail_cond_iid}, $\displaystyle\lim_{t\to \infty}\eps^{\rho}\pr{t}=0$.
Using the function $\eps^{\rho}\pr{\cdot}$ to estimate the tail probability via the elementary inequality
$\PP\pr{\abs{X_0}>t}\leq t^{-q\pr{p,\gamma}}\eps^{\rho}\pr{t}$ and the fact that $\eps^{\rho}\pr{\cdot}$ is non-increasing, one has that 
\begin{align*}
 I_{1}^{\rho}\pr{J,N,\eta}&\leq 2^N\sum_{j=J}^{N} \int_{0}^{1}u^{p-1}
 \pr{ \eta u^{1/\pr{1-q\pr{p,\gamma}}}
2^{N/2}2^{j\pr{1/p-\gamma}}}^{-q\pr{p,\gamma}}\eps^{\rho}\pr{
\eta u^{1/\pr{1-q\pr{p,\gamma}}}
2^{N/2}2^{j\pr{1/p-\gamma}}}\mathrm{d}u\\
&\leq 2^N\sum_{j=1}^{N} \int_{0}^{1}u^{p-1}
 \pr{ \eta u^{1/\pr{1-q\pr{p,\gamma}}}
2^{N/2}2^{j\pr{1/p-\gamma}}}^{-q\pr{p,\gamma}}\eps^{\rho}\pr{
\eta u^{1/\pr{1-q\pr{p,\gamma}}}
2^{N/q\pr{p,\gamma}} }\mathrm{d}u\\
&\leq  \eta^{-q\pr{p,\gamma}/\pr{q\pr{p,\gamma}-1}}\int_{0}^{1}u^{p-1+q\pr{p,\gamma}/\pr{ q\pr{p,\gamma}-1}}\eps^{\rho}\pr{
\eta u^{1/\pr{1-q\pr{p,\gamma}}}
2^{N/q\pr{p,\gamma}} }\mathrm{d}u.
\end{align*}
Hence, by the dominated convergence theorem and the fact that $\int_0^1 u^{p-1+q\pr{p,\gamma}/\pr{ q\pr{p,\gamma}-1}}\mathrm{d}u$ is finite, we claim that 
\begin{equation}\label{eq:bound_I1rho_mixing}
\forall \eta>0,\quad
\lim_{J\to\infty}\limsup_{N\to\infty} I_{1}^{\rho}\pr{J,N,\eta}=0.
\end{equation}
For the term $I_{2}(J,N)$, letting $C_\rho:=\exp\pr{\sum_{i=1}^\infty
\rho_{\mathbb{X}}\pr{2^{i}}}$ (which is finite by assumption \eqref{eq:series_coeff_rho_mixing}), we obtain 
\begin{align*}
I_{2}^{\rho}\pr{J,N} &\leq KC_\rho\left\|X_{0}\right\|_{2}^{q} \sum_{j=J}^{N}2^{j}2^{-jq\pr{1/p-\gamma}}2^{-Nq/2}  2^{(N-j)q/2}\int_{0}^{1}u^{p-1}u^{-q} \mathrm{d}u \\
& \leq KC_\rho  \left\|X_{0}\right\|_{2}^{q} \sum_{j=J}^{\infty}2^{j\frac{q\pr{p,\gamma}-q}{q\pr{p,\gamma}}} \int_{0}^{1}u^{p-q-1}\mathrm{d}u.
\end{align*}
Since $q\pr{p,\gamma}<q<p$, we infer that
\begin{equation}\label{eq:bound_I2_rho_mixing}
\lim_{J\to\infty}\limsup_{N\to\infty}I_{2}^{\rho}\pr{J,N}=0.
\end{equation}
It remains to consider $I_{3}^{\rho}\pr{J,N,\eta}$. We first bound the $q$-th moment of the truncated version of $X_0$:
\begin{align*}
\left\|X_{0}\mathbf{1}_{|X_{0}|\leq A}\right\|_{q}^{q}&
=  q\int_{0}^{\infty}t^{q-1} \mathbb{P}\left(\abs{X_0}\mathbf{1}_{|X_{0}|\leq A}> t \right)
\mathrm{d}t \\
&\leq q\int_0^At^{q-1} \mathbb{P}\left(\abs{X_0} > t \right)\mathrm{d}t \\
&\leq  qM \int_0^At^{q-q\pr{p,\gamma}-1}  \mathrm{d}t,
\end{align*}
where $M$ is defined by \eqref{eq:def_M_weal_l_q_norm},
hence
\begin{equation*}
\left\|X_{0}\mathbf{1}_{|X_{0}|\leq A}\right\|_{q}^{q}\leq \frac{Mq}{q-q\pr{p,\gamma}}A^{q-q\pr{p,\gamma}}.
\end{equation*}
Then, by letting $\pi_k:=\exp\left(K\sum_{i=1}^{k}\left(\rho_{\mathbb{X} }\pr{2^i}  \right)^{2/q} \right)$ and $K'_{p,q,\gamma}:=K Mq/\pr{q-q\pr{p,\gamma}}$, we derive the following inequality:
\begin{align*}
I_{3}^{\rho}\pr{J,N,\eta}&\leq K'_{p,q,\gamma}\sum_{j=J}^{N}2^{j}2^{N-j}\pi_{N-j}\int_0^1 u^{p-1}
\prbig{u2^{N/2}2^{j\pr{1/p-\gamma}}}^{-q} \prbig{\eta u^{\frac 1{1-q\pr{p,\gamma}}} 2^{N/2}2^{j\pr{1/p-\gamma}}}^{q-q\pr{p,\gamma}}\mathrm{d}u\\
&=\eta^{q-q\pr{p,\gamma}}2^{N\pr{1-q\pr{p,\gamma}/2}}
 K'_{p,q,\gamma}\sum_{j=J}^{N}
 \pi_{N-j}2^{-j\pr{1-q\pr{p,\gamma}/2}}\int_0^1 u^{p- \frac{q\pr{p,\gamma}\pr{q-1}}{q\pr{p,\gamma}-1}-1}
\mathrm{d}u.
\end{align*}
Since $q<p$, we have $p>\pr{q-1} q\pr{p,\gamma}/\pr{q\pr{p,\gamma}-1}$, therefore, the integral over $u$ is finite and it will be incorporated into a constant $K''_{p,q,\gamma}$. Bounding the sum over $j\in \ens{J,\dots,N}$ by that over
 $\ens{1,\dots,N}$ and then making the change of index $\ell=N-j$,
it follows that
\begin{equation*}
I_{3}^{\rho}\pr{J,N,\eta}\leq \eta^{q-q\pr{p,\gamma}}
 K''_{p,q,\gamma}\sum_{\ell=0}^{N-1}
 \pi_{\ell}2^{\ell \pr{1-q\pr{p,\gamma}/2}}.
\end{equation*}
Keeping in mind that $q\pr{p,\gamma}>2$, we infer by the ratio test that the series $s:=\sum_{\ell=0}^{\infty}
 \pi_{\ell}2^{\ell \pr{1-q\pr{p,\gamma}/2}}$ is finite. Hence, we deduce
 \begin{equation}\label{eq:bound_I3_rho_mixing}
 \limsup_{J\to\infty}\limsup_{N\to\infty}I_{3}^{\rho}\pr{J,N,\eta}\leq \eta^{q-q\pr{p,\gamma}}
  K''_{p,q,\gamma}s.
 \end{equation}
 The combination of \eqref{eq:decomposition_proba_rho_mixing} with
 \eqref{eq:bound_I1rho_mixing}, \eqref{eq:bound_I2_rho_mixing} and
 \eqref{eq:bound_I3_rho_mixing} shows that for each positive $\eta$,
 \begin{equation*}
  \lim_{J\to\infty}\limsup_{N\to\infty}\sum_{j=J}^N
 2^j\int_0^1 u^{p-1}\PP\pr{\max_{1\leq k\leq 2^{N-j}}\abs{S_{k,X_0}}>
 2^{N/2}2^{j\pr{1/p -\gamma}}  u }\mathrm{d}u
 \leq  \eta^{q-q\pr{p,\gamma}}
  K''_{p,q,\gamma}s.
 \end{equation*}
Since $\eta$ can be taken arbitrarily small and $q>q\pr{p,\gamma}$, the proof of Theorem~\ref{thm:rho_mixing} is complete.

\subsection{Proof of Theorem~\ref{thm:tau_dependent}}
The condition \eqref{eq:main-condition for tau} implies that the integral $\int_{0}^{1} Q_{X_{0}}^2(u) \cdot (\tau/2)^{-1} \circ G_{X_{0}}^{-1}(u) \, \mathrm{d}u$ is finite, which  
shows that $\mathbb{X}$ satisfies condition (5.5) from \cite{MR2199291}. For each $k\geq 1$, $\theta(k)\leq \tau(k)$ (see  \cite{MR2105738} for a definition of the coefficients $\theta\pr{k}$ and \cite[Lemma~7]{MR2105738} for a proof of this inequality). Therefore, 
 \cite[Corollary 1]{zbMATH02233768} demonstrates that   $\mathbb{X}$ satisfies the projective criterion established in \cite{MR1743095}, which in turn ensures the convergence of the finite-dimensional distributions.

It remains to establish the tightness of $\pr{W_{n,\mathbb{X}}\pr{\cdot}}_{n\geq 1}$ in 
$B_{p,\gamma}\pr{[0,1]}$. To this end, we apply Proposition~\ref{prop:deviation_tau_dep} with $\ell:=N-j$, $x:=u2^{N/2}2^{j(1/p-\gamma)}$ and $r\in \pr{q\pr{p,\gamma},p}$ and get
\begin{equation}\label{eq:decomposition_proba_tau_dep}
\sum_{j=J}^{N}2^{j}\int_{0}^{1}u^{p-1} \mathbb{P}\prbig{\max_{1\leq i \leq 2^{N-j}} \left|S_{i,X_{0}} \right|> 5u 2^{N/2}2^{j(1/p-\gamma)}}\mathrm{d}u
\leq I_{1}^{\tau}(J,N)+I_{2}^{\tau}(J,N),
\end{equation}
where
\begin{align*}
I_1^{\tau}\pr{J,N}&:=4\sum_{j=J}^{N}2^{j}\int_{0}^{1}u^{p-1}\prbig{1+\frac{(u2^{N/2}2^{j(1/p-\gamma)})^{2}}{rs_{2^{N-j}}^{2}}}^{-r/2}\mathrm{d}u, \\
I_2^{\tau}\pr{J,N}&:=\sum_{j=J}^{N}2^{j}\int_{0}^{1}u^{p-1} \frac{4\cdot 2^{N-j}}{2^{N/2}2^{j(1/p-\gamma)}u} \int_{0}^{1} Q_{X_0}(v)\mathbf{1}_{R_{\pr{\tau}}(v)\geq u2^{N/2}2^{j(1/p-\gamma)}/r}\mathrm{d}u\mathrm{d}v.
\end{align*}
Consider $I_{1}^{\tau}(J,N)$. We first have the upper bound
$$ s_{2^{N-j}}^{2}(X_{0}) \leq 4\cdot 2^{N-j}  \int_{0}^{\Vert X_{0} \Vert_1} (\tau/2)^{-1}(v) Q_{X_{0}} \circ G_{X_{0}}(v) \, \mathrm{d}v$$
from which, due to the finiteness of the integral
$\int_{0}^{\Vert X_{0} \Vert_1} (\tau/2)^{-1}(v) Q_{X_{0}} \circ G_{X_{0}}(v) \, \mathrm{d}v$, we conclude
$$ s_{2^{N-j}}^{r}(X_{0}) \leq C_{r} 2^{(N-j)r/2}.$$
This yields
\begin{align*}
I_{1}^{\tau}(J,N) & \leq C_{r}\sum_{j=J}^{N}2^{j}2^{(N-j)r/2}2^{-r(N/2+j(1/p-\gamma))}\int_{0}^{1}u^{p-r-1}\mathrm{d}u\\
& \leq \frac{C_{r}}{ p-r} \sum_{j=J}^{\infty}2^{j\pr{1-r/q\pr{p,\gamma}}}.
\end{align*}
Since $r>q(p,\gamma)$,  the series $\sum_{j=0}^{\infty}2^{j\pr{1-r/q\pr{p,\gamma}}}$ is finite, hence
\begin{equation}\label{eq:conv_H1_tau_dep}
    \lim_{J\to \infty}\limsup_{N\to \infty}I_{1}^{\tau}(J,N)=0.
\end{equation}
Consider $I_{2}^{\tau}(J,N)$. Define the function $\eps^{\tau}\colon\R_+\to\R$ by 
\begin{equation*}
\eps^{\tau}\pr{t}:=\sup_{x\geq t}x^{q(p,\gamma)-1} \int_{0}^{1} Q_{X_{0}}(u) \mathbf{1}_{R_{\pr{\tau}}\pr{u}>x} \,\mathrm{d}u.
\end{equation*}
Note that by \eqref{eq:main-condition for tau}, $\displaystyle\lim_{t\to \infty}\eps^{\tau}\pr{t}=0$.
We bound the integral over $v$ using $\eps^{\tau}\pr{\cdot}$, and we exploit the non-increasing property of $\eps^{\tau}\pr{\cdot}$ in order to bound its argument independently of $j$, we infer that
\begin{align*}
I_{2}^{\tau}(J,N)& \leq 2^{2+\frac N2}\sum_{j=1}^N\int_0^1
u^{p-2}2^{j\pr{\gamma-1/p}}\eps^\tau\prbig{\frac{u2^{N/2}2^{j\pr{1/p-\gamma}}}r}
\prbig{\frac{u2^{N/2}2^{j\pr{1/p-\gamma}}}{r}}^{1-q\pr{p,\gamma}}\mathrm{d}u\\
&\leq 2^{2+N\pr{1-\frac{q\pr{p,\gamma}}2}}\sum_{j=1}^N
2^{j\big(\gamma-1/p\big)\pr{1+q\pr{p,\gamma}-1}}
\int_0^1u^{p-q\pr{p,\gamma}-1}\eps^\tau\prbig{\frac{u2^{N/q\pr{p,\gamma}}}r}\mathrm{d}u
.
\end{align*}
Rearranging the exponents and making the change of index $k=N-j$, we obtain
\begin{equation*}
I_{2}^{\tau}(J,N)\leq C_{p,\gamma}\int_0^1u^{p-q\pr{p,\gamma}-1}\eps^\tau\prbig{\frac{u2^{N/q\pr{p,\gamma}}}r}\mathrm{d}u.
\end{equation*}
By the dominated convergence theorem,
\begin{equation}\label{eq:conv_H2_tau_dep}
    \lim_{J\to \infty}\limsup_{N\to \infty}I_{2}^{\tau}(J,N)=0.
\end{equation}
We conclude the proof of Theorem~\ref{thm:tau_dependent} by combining \eqref{eq:decomposition_proba_tau_dep} with  \eqref{eq:conv_H1_tau_dep} and \eqref{eq:conv_H2_tau_dep}.
\begin{appendices}
\section{Appendix: deviation inequalities}\label{sec:deviation_inequalities}
In the proof under the Hannan type condition \eqref{eq:Hannan}, one needs to control the contribution in the partial sum process of the tail sum of projectors. These ones are martingale differences. To do so, we need a deviation inequality for martingale difference sequences. This corresponds to \cite[Theorem~1.7]{MR4046858} stated here in the $\R$-valued case, without maxima of partial sum and that where the increments are weighted identically distributed 
martingale differences.
\begin{Proposition}[{\cite[Theorem~1.7]{MR4046858}}]
\label{devationinequality}
For each $q>2$, there exists a positive constant $C_q$, such that for each $\R$-valued stationary martingale difference sequence
$(D_i,\mathcal{F}_i)_{i\geq 0}$, $k\geq 1$, $c_1,\dots,c_k\in [-1,1]$ and $x>0$, the following inequality holds:
\begin{equation*}
    \PP\prbig{ \Big\lvert \sum_{i=1}^kc_iD_i \Big\rvert>\sqrt{k}x }
\leq C_q \int_0^\infty\PP\pr{\abs{D_0}>xv}\min\ens{v^{q-1},v}\mathrm{d}v.
\end{equation*}

\end{Proposition}

Deviation inequalities are a very convenient tool to check the tightness
criterion  mentioned in Proposition~\ref{prop:tightness_criterion}.
We collect such inequalities obtained in the literature. Let us start by
$\alpha$-dependent sequences. We refer the reader to Definitions~\ref{def:alpha_dep1}, \ref{def:alpha_dep2} and \ref{def:alpha_dep3}.
\begin{Proposition}[{\cite[Proposition~5.1]{MR3583992}}]\label{prop:deviation_alpha_dep}
Let $\mathbb{Y}=\pr{Y_i}_{i\in\Z}$, where $Y_i=Y_0\circ \mathbb{T}^i$. Denote by $\PP_{Y_0}$ the distribution of $Y_0$ and let $Q$ be a quantile function in $\el_1$. Let $X_i=f\pr{Y_i}-\E{f\pr{Y_i}}$, where $f$ belongs to the class $\widetilde{\Fca}\pr{Q,\PP_{Y_0}}$. Given a positive integer $\ell$, define
\begin{equation*}
    R_\ell\pr{v}=\min\ens{\alpha^{-1}_{2,\mathbb{Y} }\pr{v},2^\ell }Q\pr{v},\quad v\in [0,1],
\end{equation*}
and
\begin{equation*}
L_\ell\pr{x}=R_n^{-1}\pr{x}=\inf\ens{v\in [0,1]:R_{\ell}\pr{v}\leq x},\quad x\geq 0.
\end{equation*}
Let $S_{k,X_0}=\sum_{i=1}^k X_i$. Then for any $x>0$, $r>2$, $\beta\in \pr{r-2,r}$ and $\ell\geq 1$, the following inequality takes  place:
\begin{multline}\label{eq:deviation_alpha_dep}
    \PP\prbig{\max_{1\leq k\leq 2^\ell}\abs{S_{k,X_0}}  >x }
    \leq C_{r,\beta}\frac{2^{r\ell/2}}{x^r}
    \prbig{\sum_{k=0}^{2^\ell-1} \int_0^{\alpha_{1,\mathbf{Y}}\pr{k}  } Q^{2}\pr{v}\mathrm{d}v}^{r/2}+ C_{r,\beta}\frac{2^\ell}{x}
    \int_{0}^{L_{\ell}\pr{x}}Q\pr{v}\mathrm{d}v\\
    +C_{r,\beta}\frac{2^\ell}{x^{1+\beta/2}}
    \int_{0}^{L_{\ell}\pr{x}}R_{2^\ell}^{\beta/2}\pr{v}Q\pr{v}\mathrm{d}v +C_{r,\beta}\frac{2^\ell}{x^{1+r/2}}
     \int_{L_{\ell}\pr{x}}^{1}R_{\ell}^{r/2}\pr{v}Q\pr{v}\mathrm{d}v.
\end{multline}
\end{Proposition}
For $\rho$-mixing sequences, the next inequality is due to Shao.
\begin{Proposition}[{\cite[Theorem~1.2]{MR1334179}}]\label{prop:deviation_rho_mixing}
Let $\mathbb{X}=\pr{X_{i}}_{i\in\Z}$ be a strictly stationary
centered sequence and $q\geq 2$. Then there exists a constant $K$ depending only on $\pr{\rho_{\mathbb{X}}\pr{k}}_{k\geq 1}$ and $q$ such that for each $x>0$, $\ell\geq 1$ and each $A$ such that $2^{\ell+1}\norm{ X_0\mathbf{1}_{\abs{X_0}\geq A}}_1\leq x$,
\begin{multline}\label{eq:deviation_rho_mixing}
\PP\prbig{\max_{1\leq k\leq 2^\ell}\abs{S_{k,X_0}}>x}\leq
2^\ell\PP\pr{\abs{X_0}\geq A}+K x^{-q}
 2^{\ell q/2}\exp\prbig{K\sum_{i=1}^\ell\rho_{\mathbb{X}}\pr{2^i}}  \norm{X_0}_2^q \\+K x^{-q}2^\ell\exp\prbig{\sum_{i=1}^\ell \pr{\rho_{\mathbb{X}}\pr{2^i}}^{2/q}}\norm{X_0\mathbf{1}_{\abs{X_0}\leq A}}_{q}^{q}.
\end{multline}
\end{Proposition}
We now state deviation inequalities for $\tau$-dependent sequence.
For an integrable random variable  $Y$, we denote by $G_{Y}$ the generalized inverse of $x \mapsto \int_{0}^{x} Q_{Y}(u)\, \mathrm{d}u$.

The following Fuk–Nagaev inequality was established by Dedecker and Prieur.
\begin{Proposition}[{\cite[Theorem~2]{MR2199291}}]
\label{prop:deviation_tau_dep}
Let $\mathbb{X}=(X_{i})_{i \in \mathbb{Z}}$ be a strictly stationary sequence of centered and square integrable $\R$-valued random variables. Let $R_{\pr{\tau}} := ((\tau/2)^{-1} \circ G_{X_{0}}^{-1}) Q_{X_{0}}$ and $S_{\pr{\tau}}:=R^{-1}$. Then for any $x>0$, any integer $\ell \geqslant 0$ and any $r \geqslant 1$,
$$
\PP\prbig{\max_{1\leq k\leqslant 2^\ell} |S_{k,X_{0}}| \geq 5x} \leq 4 \prbig{1+\frac{x^{2}}{r s_{2^\ell}^2\pr{\mathbb{X}}} }^{-r/2}+\frac{2^{\ell+2}}{x} \int_{0}^{S_{\pr{\tau}}(x/r)} Q_{X_{0}}(v) \, \mathrm{d}v,
$$
where 
$$
s_{2^\ell}^2\pr{\mathbb{X}}:=
\sum_{i=1}^{2^\ell} \sum_{j=1}^{2^\ell} \left| \operatorname{Cov}(X_{i}, X_{j})\right|
\leq 2^{\ell+2} \int_{0}^{\Vert X_{0} \Vert_1}(\tau/2)^{-1}(v) Q_{X_{0}} \circ G_{X_{0}}(v) \, \mathrm{d}v.
$$
\end{Proposition}
\end{appendices}

\subsection*{Acknowledgments} 
This work was finalized during the second author's visit to the University of 
Strasbourg. The author is grateful to IMU-CDC, the Simons Foundation and to University of 
Strasbourg for the hospitality.
\def\cprime{$'$}
\providecommand{\bysame}{\leavevmode\hbox to3em{\hrulefill}\thinspace}
\providecommand{\MR}{\relax\ifhmode\unskip\space\fi MR }
\providecommand{\MRhref}[2]{%
  \href{http://www.ams.org/mathscinet-getitem?mr=#1}{#2}
}

\end{document}